\documentclass{siamart0516}
\usepackage{tikz}
\usepackage{arydshln}
\usepackage{comment}
\usepackage{esint}
\usepackage{stmaryrd}
\usepackage{mathtools}
\usetikzlibrary{shapes}

\def\NF{{\mathcal F}}
\def\WS{{\mathcal W}}
\def\meq{{M_\text{eqn}}}
\def\mdim{{M_\text{dim}}}
\def\mdeg{{M_\text{deg}}}
\def\melems{{M_\text{elem}}}
\def\mbasis{M_{\text{basis}}}
\def\mpred{{M_\text{P}}}
\def\mcorr{{M_\text{C}}}
\def\hskipper{\hspace{-1.1mm}}
\def\Lone{\mat{L^{\text{   \hspace{-2.5mm}  1 \hskipper 0 \hskipper 1 \hskipper 0}}}}
\def\Ltwo{\mat{L^{\text{   \hspace{-2.5mm}  1 \hskipper 1 \hskipper 1 \hskipper 0}}}}
\def\Lthree{\mat{L^{\text{ \hspace{-2.2mm}  0 \hskipper 1 \hskipper 1 \hskipper 0}}}}
\def\Lfour{\mat{L^{\text{  \hspace{-2.5mm}  1 \hskipper 0 \hskipper 1 \hskipper 1}}}}
\def\Lfive{\mat{L^{\text{  \hspace{-2.5mm}  1 \hskipper 1 \hskipper 1 \hskipper 1}}}}
\def\Lsix{\mat{L^{\text{   \hspace{-2.2mm}  0 \hskipper 1 \hskipper 1 \hskipper 1}}}}
\def\Lseven{\mat{L^{\text{ \hspace{-2.5mm}  1 \hskipper 0 \hskipper 0 \hskipper 1}}}}
\def\Leight{\mat{L^{\text{ \hspace{-2.5mm}  1 \hskipper 1 \hskipper 0 \hskipper 1}}}}
\def\Lnine{\mat{L^{\text{  \hspace{-2.2mm}  0 \hskipper 1 \hskipper 0 \hskipper 1}}}}

\def\beq{\begin{equation}}
\def\eeq{\end{equation}}

\def\bga{\begin{gather}}
\def\ega{\end{gather}}

\def\bal{\begin{align}}
\def\eal{\end{align}}

\def\Tm{{\mathcal T}}

\newcommand{\reals}{\mathbb R}

\newcommand{\Wp}{ \vec{W}^{n+1/2}_{i+1} }
\newcommand{\Wo}{ \vec{W}^{n+1/2}_{i} }
\newcommand{\Wm}{ \vec{W}^{n+1/2}_{i-1} }
\newcommand{\Wppast}{ \vec{W}^{n}_{i+1} }
\newcommand{\Wopast}{ \vec{W}^{n}_{i} }
\newcommand{\Wmpast}{ \vec{W}^{n}_{i-1} }
\newcommand{\psieast}{ \vec{\Psi}_{|_{\xi=1}} }
\newcommand{\psiwest}{ \vec{\Psi}_{|_{\xi=-1}} }
\newcommand{\psifutr}{ \vec{\Psi}_{|_{\tau=1}} }
\newcommand{\psipast}{ \vec{\Psi}_{|_{\tau=-1}} }
\newcommand{\psitau}{ \vec{\Psi}_{|_\tau } }
\newcommand{\psixii}{ \vec{\Psi}_{|_\xi } }

\newcommand{\bunderline}[1]{\underline{#1}}
\renewcommand{\vec}[1]{{\bunderline{#1}}}
\newcommand{\vect}[1]{{\bunderline{#1}}}
\newcommand{\mat}[1]{{\bunderline{\bunderline{#1}}}}

\newcommand{\abs}[1]{\left| #1 \right|}

\newcommand{\paren}[1]{\left( #1 \right)}

\newcommand{\vectthree}[3]{\left(\begin{array}{c} #1\\ #2 \\ #3 \end{array} \right)}

\newcommand{\half}{\frac{1}{2}}

\newcommand{\grad}{\vec{\nabla}}
\renewcommand{\div}{\vec{\nabla} \cdot}

\newsiamthm{remark}{Remark}

\makeatletter
\renewcommand*\env@matrix[1][\arraystretch]{%
  \edef\arraystretch{#1}%
  \hskip -\arraycolsep
  \let\@ifnextchar\new@ifnextchar
  \array{*\c@MaxMatrixCols c}}
\makeatother

\usepackage{lipsum}
\usepackage{amsfonts}
\usepackage{amssymb}
\usepackage{graphicx}
\usepackage{epstopdf}
\usepackage{algorithmic}
\ifpdf
  \DeclareGraphicsExtensions{.eps,.pdf,.png,.jpg}
\else
  \DeclareGraphicsExtensions{.eps}
\fi

\newcommand{\TheTitle}{The Regionally-Implicit Discontinuous Galerkin Method:
Improving the Stability of DG-FEM}
\newcommand{\TheAuthors}{Pierson T. Guthrey and James A. Rossmanith}

\headers{Regionally-Implicit Discontinuous Galerkin Method}{\TheAuthors}

\title{{\TheTitle}\thanks{Submitted to the editors on {\today}. \funding{This work was funded in part by NSF Grant DMS--1620128.}}}

\author{
  Pierson T. Guthrey\thanks{Michigan State University, Department of
  Computational Mathematics, Science and Engineering, 428 S. Shaw Lane, East Lansing, Michigan 48824, USA
    (\email{piersonguthrey@gmail.com}).}
  \and
  James A. Rossmanith\thanks{Iowa State University, Department of
  Mathematics, 411 Morrill Road, Ames, Iowa 50011, USA
    (\email{rossmani@iastate.edu}).}
}

\usepackage{amsopn}

\ifpdf
\hypersetup{
  pdftitle={\TheTitle},
  pdfauthor={\TheAuthors}
}
\fi

\begin{document}

\maketitle

\begin{abstract}
Discontinuous Galerkin (DG) methods for hyperbolic partial differential equations (PDEs) with explicit time-stepping schemes, such as strong stability-preserving Runge-Kutta (SSP-RK), suffer from time-step restrictions that are significantly worse than what a simple Courant-Friedrichs-Lewy (CFL) argument requires. In particular, the maximum stable time-step scales inversely with the highest degree in the DG polynomial approximation space and becomes progressively smaller with each added spatial dimension. In this work we introduce a novel approach that we have dubbed the regionally implicit discontinuous Galerkin (RIDG) method  to overcome these small time-step restrictions. The RIDG method is based on an extension of the Lax-Wendroff DG (LxW-DG) method, which previously had been shown to be equivalent  to a predictor-corrector approach, where the predictor is a locally implicit spacetime method (i.e., the predictor is something like a block-Jacobi update for a fully implicit spacetime DG method). The corrector is an explicit method that uses the spacetime reconstructed solution from the predictor step. In this work we modify the predictor to include not just local information, but also neighboring information. With this modification we show that the stability is greatly enhanced; in particular, we show that we are able to remove the polynomial degree dependence of the maximum time-step and show how this extends to multiple spatial dimensions. A semi-analytic von Neumann analysis is presented to theoretically justify the stability claims. Convergence and efficiency studies for linear and nonlinear problems in multiple dimensions are accomplished
using a {\sc matlab} code that can be freely downloaded.
\end{abstract}

\begin{keywords}
  discontinuous Galerkin, hyperbolic conservation laws, Courant-Friedrichs-Lewy condition, time-setpping, numerical stability
\end{keywords}

\begin{AMS}
  65M12, 65M60, 35L03
\end{AMS}

\section{Introduction}
\label{sec:intro}
Hyperbolic conservation laws model phenomena characterized by waves propagating at
finite speeds; examples include the shallow water (gravity waves), compressible Euler (sound waves), Maxwell (light waves), magnetohydrodynamic (magneto-acoustic and Alfv\'en waves), and Einstein (gravitational waves) equations. 
In recent years, the discontinuous Galerkin (DG) finite element method (FEM) has become a standard approach for solving hyperbolic conservation laws alongside other methods such as weighted essentially non-oscillatory (WENO) schemes (e.g., see Shu \cite{article:ShuWENO2009}) and 
various finite volume methods (e.g., see LeVeque \cite{book:Le02}). The DG method
was first introduced by Reed and Hill \cite{article:ReedHill73} for neutron transport, and then fully developed for time-dependent hyperbolic conservation laws in a series of papers by  Cockburn,  Shu, and collaborators (see \cite{article:CoShu98} and references therein for details).
An important feature of DG methods is that they can, at least in principle, be made arbitrarily high-order in space by  increasing the polynomial order in each element; and therefore, the DG method is an example of a spectral element method (e.g., see Chapter 7.5 of Karniadakis and Sherwin \cite{book:KaSh2005}). 

If DG is only used to discretize the spatial part of the underlying PDE, it remains to also introduce a temporal discretization. Many time-stepping methods are possible, including various explicit and implicit schemes. In general, one time-step of an explicit  scheme is significantly cheaper than an implicit one; the trade-off is that 
implicit schemes usually allow for larger time-steps. In many applications involving hyperbolic conservation laws, however, it is necessary to resolve the fastest time scales, in which case explicit methods are more efficient and easier to implement than implicit ones. 

An upper bound on the largest allowable time-step for explicit schemes is provided by the Courant-Friedrichs-Lewy (CFL) condition, which requires that the domain of dependence of the numerical discretization subsumes the domain of dependence of the continuous PDE \cite{article:CFL1928}.
For example, a 1D hyperbolic PDE for which information propagates at a maximum wave speed of $\lambda_{\text{max}}$, on a uniform mesh of elements of size $h=\Delta x$, and with a time-stepping method that updates the solution on
the element $\Tm_i^h$ only using existing solution values from $\Tm_{i-1}^h$, $\Tm_{i}^h$, and $\Tm_{i+1}^h$,  
 has the following constraint on  
$\Delta t$:
\begin{equation}
\nu := \frac{\lambda_{\text{max}} \Delta t}{\Delta x} \le 1.
\end{equation}
This has a clear physical interpretation: a wave that emanates from the boundaries of element $i$ that is
traveling at the maximum speed, $\lambda_{\text{max}}$, is not allowed to propagate further than one element width.
If we wanted to allow the wave to travel more than one element width, we would need to widen the numerical stencil.

The CFL condition as described above is a {\it necessary} condition for stability (and therefore convergence), but it is not {\it sufficient}. For high-order DG methods with explicit time-stepping, a fact that is well-known in the literature is that the actual maximum linearly stable value of the CFL number, $\nu = {\lambda_{\text{max}} \Delta t}/{\Delta x}$, is significantly smaller than what the CFL condition predicts (see for example Liu et al. \cite{article:LiuShuTadmorZhang08} and Sections 4.7 and 4.8 of Hesthaven and Warburton \cite{book:HesWar2007}).

Two popular explicit time-stepping schemes for DG are strong-stability-preserving Runge-Kutta DG (SSP-RK) \cite{article:GoShu98, gottliebShuTadmor01} and Lax-Wendroff \cite{article:GasDumHinMun2011,article:QiuDumShu2005}.  SSP-RK time-steps are one-step multistage Runge-Kutta methods that can be written as convex combinations of forward Euler steps. Lax-Wendroff utilizes the Cauchy-Kovalevskaya 
\cite{article:Kovaleskaya1875} procedure to convert temporal derivatives into spatial derivatives; the name Lax-Wendroff is due to the paper  of Lax and Wendroff \cite{article:LxW1960}. 
In \cref{table:CFL_gap} we illustrate for both the SSP-RK and Lax-Wendroff DG methods the gap between the CFL condition, a necessary but not sufficient condition for stability, and the 
semi-analytically computed maximum CFL number needed for linear stability.
Shown are the methods with space and time order $k=1,2,3,4$.
The SSP-RK DG numbers are reported from Liu et al. \cite{article:LiuShuTadmorZhang08}, while the Lax-Wendroff DG numbers are from von Neumann analysis done in this paper. Note that the maximum CFL number from the CFL condition for SSP-RK DG grows with $k$ due to the fact that the number of Runge-Kutta stages grows with $k$; and therefore, the numerical domain of dependence is increased. For both sets of methods, the clear trend is that the maximum CFL numbers are much smaller than what a simple CFL domain of dependence argument would dictate.
In particular, the relationship between the maximum CFL number 
and the order of the method is roughly: $\nu_{\text{max}} \propto 1/k$.

The goal of this paper is to develop an alternative time discretization for DG that allows for a linearly stable time-step that is closer to what is predicted by the CFL condition. The starting point of this work is the
interpretation of the Lax-Wendroff DG method developed by Gassner et al. \cite{article:GasDumHinMun2011}, where it was shown that Lax-Wendroff DG can be formulated as a predictor-corrector method. The predictor is a local version of a spacetime DG method \cite{article:KlaVegVen2006,article:Sudirham2006} 
(i.e., the predictor is something like a block-Jacobi update for a fully implicit spacetime DG method), and the corrector is an explicit method that uses the spacetime reconstructed solution from the predictor step.
In this work we modify the predictor to include not just local information, but also neighboring information.
The name that we are giving to this new approach is the {\it regionally-implicit} discontinuous Galerkin (RIDG) scheme,
which contrasts with the {\it locally-implicit} (LIDG) formulation of the Lax-Wendroff DG scheme developed by
Gassner et al. \cite{article:GasDumHinMun2011}.
In this new formulation, we are able to achieve all of the following:
\begin{itemize}
\item Develop RIDG schemes for 1D, 2D, and 3D advection;
\item Show that RIDG has larger maximum CFL numbers than explicit SSP-RK and Lax-Wendroff DG;
\item Show that the maximum linearly stable CFL number is bounded below by a constant that is independent of the polynomial order;
\item Demonstrate experimentally the correct convergence rates on 1D, 2D, and 3D advection examples.
\item Demonstrate experimentally the correct convergence rates on 1D and 2D nonlinear examples.
\end{itemize}
All of the methods described in this work are written in a {\sc matlab} code that can be freely downloaded
\cite{code:ridg-code}.
 
The organization of this paper is as follows.  In \cref{sec:dg-fem-space} we briefly review how space is discretized in the discontinuous Galerkin (DG) method. In \cref{sec:one-dimension} we review the Lax-Wendroff DG scheme, then develop the one-dimensional version of the proposed regionally implicit DG (RIDG) scheme, and carry out von Neumann stability analysis for both methods. 
The generalization to multiple dimensions is done in \cref{sec:higher-dimensions}. In \cref{sec:results} we carry out numerical convergence tests to validate the new approach and
to quantify the computational efficiency of RIDG relative to the Lax-Wendroff method.
Finally, in \cref{sec:burgers} we show how to extend the method to a nonlinear scalar problem: the 1D and 2D Burgers equation. In this case we compare the efficiency and accuracy of our method against the fourth-order Runge-Kutta discontinous Galerkin  (RKDG) scheme.

\begin{table}[!t]
\begin{center}
\begin{tabular}{|c||c|c|c|c||c|c|c|c|}
\hline 
& \multicolumn{4}{|c||}{{\bf SSP-RK($k$) with $P^{k-1}$ DG}} 
& \multicolumn{4}{|c|}{{\bf $P^{k-1}$ Lax-Wendroff DG}}\\
\hline\hline
 & $k=1$ & $k=2$ & $k=3$ & $k=4$ & $k=1$ & $k=2$ & $k=3$ & $k=4$ \\
\hline
{\bf CFL cond.} & 1.00 & 2.00 & 3.00 & 4.00 & 1.000 & 1.000 & 1.000 & 1.000 \\
\hline
{\bf Neumann} & 1.00 & 0.33 & 0.13 & 0.10  & 1.000 & 0.333 & 0.171 & 0.104  \\
\hline
\end{tabular}
\caption{Shown here are the maximum CFL numbers for the SSP-RK DG and Lax-Wendroff DG methods with the same time and space order of accuracy. The line labeled ``CFL cond.'' is the upper bound of the CFL number as predicted just by looking at the domain of dependence of the numerical method. The line label ``Neumann'' is the numerically calculated maximum linearly stable CFL number. The SSP-RK DG numbers are reported from Liu et al. \cite{article:LiuShuTadmorZhang08}, while the Lax-Wendroff DG numbers are from von Neumann analysis done in this paper. For both sets of methods, the clear trend is that the maximum linearly stable CFL numbers are much smaller than what a simple CFL domain of dependence argument would dictate.
\label{table:CFL_gap}}
\end{center}
\end{table}

\section{DG-FEM spatial discretization}
\label{sec:dg-fem-space}
Consider hyperbolic conservation laws of the form
\begin{equation}
\label{eqn:conslaw}
\vec{q}_{,t} + \div \mat{F}\left( \vec{q} \right) = \vec{0},
\end{equation}
where $\vec{q}\left(t,\vec{x} \right): \reals^+ \times \reals^\mdim \mapsto \reals^\meq$ is the vector of conserved variables, $\mat{F}\left(\vec{q}\right): \reals^\meq \mapsto \reals^{\meq \times \mdim}$ is the flux function, $\mdim$ is the number of spatial dimensions, and $\meq$ is the number of conserved variables. We assume that the system is hyperbolic, which means that the flux Jacobian,
\begin{equation}
\mat{A}\left( \vec{q}; \vec{n} \right) =  \frac{\partial \left( \vec{n} \cdot  \mat{F} \right)}{\partial \vec{q}},
\end{equation}
 for all $\vec{q} \in {\mathcal S} \subset \reals^\meq$, where ${\mathcal S}$ is some physically meaningful convex subset of $\reals^\meq$,
 and for all directions, $\vec{n} \in \reals^\mdim$ such that $\| \vec{n} \| = 1$, must be diagonalizable with only
 real eigenvalues (e.g., see Chapter 18 of LeVeque \cite{book:Le02}).

Next consider discretizing system \cref{eqn:conslaw} in space via the discontinuous Galerkin (DG) method, which
was first introduced by Reed and Hill \cite{article:ReedHill73} for neutron transport, and then fully developed for time-dependent hyperbolic conservation laws in a series of papers by Bernardo Cockburn, Chi-Wang Shu, and collaborators (see \cite{article:CoShu98} and references therein for details). We define $\Omega \subset \reals^\mdim$ to be a polygonal domain with 
boundary $\partial \Omega$, and discretize $\Omega$ using a finite set of non-overlapping elements, $\Tm_i$, such
that $\cup_{i=1}^\melems \Tm_i = \Omega$, where $\melems$ is the total number of elements.
Let ${\mathbb P}\left(\mdeg, \mdim \right)$ denote the set of polynomials from $\reals^\mdim$ to $\reals$ 
with maximal polynomial degree $\mdeg$\footnote{In 1D (i.e., $\mdim=1$), this definition is unambiguous. In higher dimensions,
${\mathbb P}\left(\mdeg, \mdim \right)$ could refer the set of polynomials that have a total degree $\le \mdeg$ (we refer to this
as the ${\mathcal P}\left(\mdeg, \mdim \right)$ basis), it could refer to the set of polynomials that have degree $\le \mdeg$
in each independent variable (we refer to this as the ${\mathcal Q}\left(\mdeg, \mdim \right)$ basis), or it could be something in between.}.
On the mesh of $\melems$ elements we define the {\it broken} finite element space:
\begin{equation}
\label{eqn:broken_space}
    \WS^h := \left\{ \vec{w}^h \in \left[ L^{\infty}(\Omega) \right]^{\meq}: \,
    \vec{w}^h \bigl|_{\Tm_i} \in \left[ {\mathbb P} \left(\mdeg, \mdim \right) \right]^{\meq} \, \, \forall \Tm_i \right\},
\end{equation}
where $h$ is the grid spacing, $M_{\text{dim}}$ is the number of spatial dimensions, 
$\meq$ is the number of conserved variables, and $\mdeg$ is the 
maximal polynomial degree in the finite element representation.
The above expression means that $\vec{w} \in \WS^h$ has
$\meq$ components, each of which when restricted to some element $\Tm_i$
is a polynomial in ${\mathbb P}\left(\mdeg, \mdim \right)$, and no continuity is assumed
across element faces. 

Let $\varphi_{k}\left(\vec{x}\right)$ for $k=1,\ldots,\mbasis$ be an appropriate basis that spans 
${\mathbb P} \left(\mdeg, \mdim \right)$ over $\Tm_i$ (e.g., Legendre or Lagrange polynomials).
In order to get the DG semi-discretization, we multiply \cref{eqn:conslaw} by 
$\varphi_{k} \in {\mathbb P} \left(\mdeg, \mdim \right)$, integrate over the element $\Tm_i$, use integration-by-parts in space,
and replace the true solution, $\vec{q}$, by the following ansatz:
\begin{equation}
\vec{q}^h\left(t, \vec{x} \right) \Bigl|_{\Tm_i} = \sum_{\ell=1}^{\mbasis} \vec{Q}_{i}^{\ell}(t) \, \varphi_{\ell}\left(\vec{x}\right).
\end{equation}
All of these steps results in the following semi-discrete system:
\begin{equation}
\label{eqn:semi_discrete_dg}
\sum_{\ell=1}^{\mbasis} \left[ \int_{\Tm_i} \varphi_{k} \varphi_{\ell} \, d\vec{x} \right] \frac{d\vec{Q}^{\ell}}{dt}
= \int_{\Tm_i} \mat{F}\left( \, \vec{q}^h \right) \cdot
\grad \varphi_{k} \, d\vec{x} - \oint_{\partial \Tm_i} \varphi_k \, \vec{{\mathcal F}}\left(\vec{q}^h_{+}, \vec{q}^h_{-};
\vec{n} \right) \, d\vec{s},
\end{equation}
where $\vec{n}$ is an outward-pointing normal vector to $\partial \Tm_i$, $\vec{q}^h_{+}$ and $\vec{q}^h_{-}$ are the states on either side of the boundary $\partial \Tm_i$, and $\vec{{\mathcal F}}$ is the numerical flux, which must satisfy
the following two conditions:
\begin{itemize}
  \item Consistency:  \, 
  $\vec{\NF}\left( \, \vec{q}, \, \vec{q}; \, \vec{n} \right) = \mat{F}\left( \, \vec{q} \, \right) \cdot \vec{n};
  $
  \item Conservation: \, 
  $\vec{\NF}\left(\vec{q}^h_-, \, \vec{q}^h_+; \, \vec{n} \right) = -\vec{\NF}\left(\vec{q}^h_+, \, \vec{q}^h_-; \, -\vec{n} \right)$.
\end{itemize}
Equation \cref{eqn:semi_discrete_dg} represents a large system of coupled ordinary differential equations in time.

\section{RIDG in one space dimension}
\label{sec:one-dimension}
We present in this section the proposed regionally-implicit discontinuous Galerkin (RIDG) method as applied to a 
one dimensional advection equation.
Each RIDG time-step is comprised of two key steps: a predictor and a corrector. The predictor is a truncated version of an implicit spacetime DG approximation, which is not consistent, at least by itself, with the PDE that it endeavors to approximate. The corrector is a modified forward Euler step that makes use of the predicted solution; this step restores consistency, and indeed, high-order accuracy, with the underlying PDE. 

In the subsections below we begin with a brief description of the advection equation in \cref{sec:RIDG1D_advection}. 
We then review the Lax-Wendroff (aka locally-implicit) prediction step in \cref{sec:LIDG1D_predict}, which provides the motivation for RIDG. The RIDG prediction step is developed in \cref{sec:RIDG1D_predict}. The correction
step for both predictors is detailed in 
\cref{sec:RIDG1D_correct}. Finally, we carry out semi-analytic von Neumann analysis for both schemes in
\cref{sec:RIDG1D_stability} and demonstrate the improved stability of RIDG over Lax-Wendroff DG.

\subsection{1D advection equation}
\label{sec:RIDG1D_advection}
We consider here the 1D advection equation for $(t,x) \in [0,T] \times \Omega$,
along with some appropriate set of boundary conditions:
\begin{equation}
\label{eqn:adv1d}
q_{,t} + u q_{,x} = 0.
\end{equation}
Next, we introduce a uniform Cartesian spacetime mesh with spacetime elements:
\begin{equation}
\label{eqn:spacetime_elem_1d}
{\mathcal S}^{n+1/2}_i = \left[t^{n}, t^n + \Delta t \right] \times \left[ x_i - {\Delta x}/{2}, 
x_i + {\Delta x}/{2} \right],
\end{equation}
which can ben written in local coordinates, $[\tau, \xi] \in [-1,1]^2$, where
\begin{equation}
t = t^{n+1/2} + \tau \left( {\Delta t}/{2} \right) \quad \text{and} \quad 
x = x_i + \xi \left( {\Delta x}/{2} \right).
\end{equation}
In these local coordinates the advection equation \cref{eqn:adv1d} becomes
\begin{equation}
\label{eqn:adv1d_nondim}
q_{,\tau} + \nu q_{,\xi} = 0, \quad \text{where} \quad \nu = \frac{u \Delta t}{\Delta x},
\end{equation}
and $|\nu|$ is the CFL number.

\subsection{Lax-Wendroff DG (aka LIDG) prediction step}
\label{sec:LIDG1D_predict}
We review here the prediction step for the Lax-Wendroff DG scheme as formulated by Gassner et al. \cite{article:GasDumHinMun2011}.
In order to contrast with the proposed RIDG method, we will refer to this method as the locally-implicit DG (LIDG) method. 

We fix the largest polynomial degree to $\mdeg$ in order to eventually achieve an approximation that has
an order of accuracy ${\mathcal O}\left(\Delta x^{\mdeg+1} + \Delta t^{\mdeg+1}\right)$.
At the old time, $t=t^n$, we are given the following approximate solution on each space element,
$\Tm_i = \left[x_i - \Delta x/2, \, x_i + \Delta x/2\right]$:
\begin{equation}
\label{eqn:old_ansatz}
q(t^{n},x) \Bigl|_{{\mathcal T}_i} \approx q^{n}_{i} := 
\vec{\Phi}^T \vec{Q}^{n}_i,
\end{equation}
where $\vec{Q}^{n}_i \in \reals^{\mcorr}$, $\vec{\Phi} \in \reals^{\mcorr}$, $\mcorr:=\mdeg+1$, and
\begin{equation}
\label{eqn:phi_basis_1d}
\vec{\Phi} = \left( 1, \, \sqrt{3} \xi, \, \frac{\sqrt{5}}{2} \left( 3 \xi^2 - 1 \right), \, \cdots \right),
\quad \text{s.t.} \quad \half \int_{-1}^{1} \vec{\Phi} \, \vec{\Phi}^T \, d\xi = \mat{\mathbb I} \in \reals^{\mcorr\times\mcorr},
\end{equation}
are the orthonormal space Legendre polynomials.

In order to compute a predicted solution on each spacetime element \cref{eqn:spacetime_elem_1d}
we make the following ansatz:
\begin{equation}
\label{eqn:pred_ansatz}
q(t,x) \Bigl|_{{\mathcal S}^{n+1/2}_i} \approx w^{n+1/2}_{i} := 
\vec{\Psi}^T \vec{W}^{n+1/2}_i,
\end{equation}
where $\vec{W}^{n+1/2}_i \in \reals^{\mpred}$, $\vec{\Psi} \in \reals^{\mpred}$, $\mpred := (\mdeg+1)(\mdeg+2)/2$, and
\begin{equation}
\label{eqn:psi_basis_1d}
\vec{\Psi} = \left( 1, \, \sqrt{3} \tau, \, \sqrt{3} \xi, \, \cdots \right), \quad
\text{s.t.} \quad \frac{1}{4} \int_{-1}^{1} \int_{-1}^{1} \vec{\Psi} \, \vec{\Psi}^T \, d\tau \, d\xi  = \mat{\mathbb I} \in \reals^{\mpred\times\mpred},
\end{equation}
are the spacetime Legendre basis functions.
Next, we pre-multiply \cref{eqn:adv1d_nondim} by $\vec{\Psi}$ and integrate over ${\mathcal S}^{n+1/2}_i$ to obtain:
\begin{equation}
\label{eqn:integrate_1d_adv}
\frac{1}{4} \int_{-1}^{1} \int_{-1}^{1} \vec{\Psi} \, \left[ q_{,\tau} +\nu q_{,\xi} \right]
\, d\tau \, d\xi = \vec{0}.
\end{equation}
We then replace the exact solution, $q$, by \cref{eqn:pred_ansatz}.
We integrate-by-parts in time, first forwards, then backwards, which introduces a jump term at the old time $t=t^n$.
No integration-by-parts is done in space -- this is what gives the local nature of the predictor step. 
All of this results in the following equation:
\begin{equation}
\begin{split}
\iint & \vec{\Psi} \, \left[ \vec{\Psi}_{,\tau} + \nu
 \vec{\Psi}_{,\xi} \right]^T \, \vec{W}^{n+1/2}_i \, d\tau \, d\xi 
+ \int  \vec{\Psi}_{|_{\tau=-1}}  \left[ \vec{\Psi}_{|_{\tau=-1}}^T \vec{W}^{n+1/2}_i
 - \vec{\Phi}^T \vec{Q}^{n}_i \right] \, d\xi = \vec{0},
\end{split}
\end{equation}
where all 1D integrals are over $[-1,1]$, which can be written as
\begin{gather}
\label{eqn:prediction_soln}
 \mat{L^0} \, \vec{W}^{n+1/2}_i = \mat{T} \, \vec{Q}^{n}_i, \\
\label{eqn:predictedA}
\mat{L^0} = \frac{1}{4} \int_{-1}^{1} \int_{-1}^{1} \vec{\Psi} \, \left[ \vec{\Psi}_{,\tau} + \nu \vec{\Psi}_{,\xi} \right]^T  \, d\tau \, d\xi 
+ \frac{1}{4}  \int_{-1}^{1} \vec{\Psi}_{|_{\tau=-1}} \, \vec{\Psi}_{|_{\tau=-1}}^T \, d\xi \in \reals^{\mpred \times \mpred}, \\
\label{eqn:predictedB}
\mat{T} = \frac{1}{4}  \int_{-1}^{1} \vec{\Psi}_{|_{\tau=-1}} \, \vec{\Phi}^T  \, d\xi \in \reals^{\mpred \times \mcorr}.
\end{gather}
As is evident from the formulas above, the predicted spacetime solution as encoded in the coefficients
$\vec{W}^{n+1/2}_i$ is entirely local -- the values only depend on the old values from the same element:
$\vec{Q}^{n}_i$. Therefore, we refer to this prediction step as {\it locally-implicit}.

\begin{remark}
Gassner et al. \cite{article:GasDumHinMun2011} argued that the locally-implicit prediction step as presented above produces the key step in the Lax-Wendroff DG scheme \cite{article:QiuDumShu2005}. We briefly illustrate this point here.

 Lax-Wendroff  \cite{article:LxW1960} (aka the Cauchy-Kovalevskaya \cite{article:Kovaleskaya1875} procedure)
   begins with a Taylor series in time. All time derivatives are then replaced by spatial derivatives using the underlying PDE -- in this case \cref{eqn:adv1d}. 
For example, if we kept all the time derivatives up to the third derivative we would get:
\begin{equation}
\begin{split}
q^{n+1} &\approx  q^n + 2 q^n_{,\tau} + 2 q^n_{,\tau,\tau} + \frac{4}{3} q^n_{,\tau,\tau,\tau}
 = q^n - 2 \nu q^n_{,\xi} + 2 \nu^2 q^n_{,\xi,\xi} - \frac{4}{3} \nu^3 q^n_{,\xi,\xi,\xi}
 \\ &= q^n - 2 \nu \left\{ q^n -  \nu q^n_{,\xi} + \frac{2}{3} \nu^2 q^n_{,\xi,\xi} \right\}_{,\xi}
 = q^n - 2 \nu {\mathcal G}^{\text{LxW}}_{,\xi}.
\end{split}
\end{equation}
Using the third-order DG approximation,
\begin{equation}
\label{eqn:third_order_dg}
q^n := \varphi_1  Q_1^n + \varphi_2  Q_2^n + \varphi_3  Q_3^n = 
Q_1^n + \sqrt{3} \xi Q_2^n + \frac{\sqrt{5}}{2} \left( 3\xi^2 -1 \right) Q_3^n,
\end{equation}
we obtain
\begin{equation}
\label{eqn:lxw_flux_lxw}
{\mathcal G}^{\text{LxW}} =  \varphi_1  \left( Q_1^n - \sqrt{3} \nu \, Q_2^n + 2 \sqrt{5} \nu^2 Q_3^n \right) + \varphi_2  \left( Q_2^n - \sqrt{15} \nu Q_3^n \right)
+ \varphi_3  Q_3^n.
\end{equation}
Alternatively, the time-averaged flux, ${\mathcal G}^{\text{LxW}}$, can be directly obtained from the locally-implicit predictor described above. We first calculate the predicted spacetime solution, $w^{n+1/2}$, via \cref{eqn:prediction_soln}, \cref{eqn:predictedA}, and \cref{eqn:predictedB}. From this we compute the time averaged flux:
\begin{equation}
\label{eqn:predict_flux_lxw}
{\mathcal G}^{\text{LxW}} = \half \int_{-1}^{1} w^{n+1/2}\left(\tau,\xi\right) d\tau 
=\half \left\{ \int_{-1}^{1} \vec{\Psi}\left(\tau,\xi\right) \, d\tau \right\}^T
\left( \mat{L^0} \right)^{-1} \mat{T} \, \vec{Q}^n.
\end{equation}
A straightforward calculation shows that \cref{eqn:predict_flux_lxw} with \cref{eqn:third_order_dg} and the corresponding 
${\mathcal P}(2,2)$ spacetime basis\footnote{Actually, for LIDG the result is the same whether we use
${\mathcal P}(\mdeg,\mdim+1)$, ${\mathcal Q}(\mdeg,\mdim+1)$, or something in between.}
is exactly the same as \cref{eqn:lxw_flux_lxw}.
\end{remark}

\begin{figure}[!t]
\centering
\begin{tikzpicture}[scale=0.68, 
circ/.style={scale=1.5, shape=circle, inner sep=2pt, draw, fill=red,node contents=},
dia/.style={scale=1.2, shape=diamond, inner sep=2pt, draw, fill=yellow,node contents=},
amp/.style = {scale=1, regular polygon, regular polygon sides=3, shape border rotate=-180,
              draw, fill=gray!30, inner sep=2pt, node contents=}]
\draw[->,line width=0.5mm, gray] (-0.5,-0.5)--(-0.5,1) node[above,black]{$t$};
\draw[->,line width=0.5mm, gray] (-0.5,-0.5)--(1,-0.5) node[right,black]{$x$};
\draw[->,line width=0.5mm, gray] (5.5,-0.5)--(5.5,1) node[above,black]{$t$};
\draw[->,line width=0.5mm, gray] (5.5,-0.5)--(7,-0.5) node[right,black]{$x$};
\draw [line width=0.8mm] (0,0) -- (4,0) -- (4,4) -- (0,4) -- cycle;
\draw [line width=0.2mm, black] node (c2) at (0, 2) [circ];
\draw [line width=0.2mm, black] node (c5) at (4, 2) [circ];
\draw [line width=0.2mm, black] node (c8) at (2, 4) [amp];
\draw [line width=0.2mm, black] node (c11) at (2, 0) [amp];
\draw [line width=0.8mm] (6,0) -- (10,0) -- (10,4) -- (6,4) -- cycle;
\draw [line width=0.8mm] (10,0) -- (14,0) -- (14,4) -- (10,4) -- cycle;
\draw [line width=0.8mm] (14,0) -- (18,0) -- (18,4) -- (14,4) -- cycle;
\draw [line width=0.2mm, black] node (r2) at (6, 2) [circ];
\draw [line width=0.2mm, black] node (r5) at (10, 2) [dia];
\draw [line width=0.2mm, black] node (r8) at (14, 2) [dia];
\draw [line width=0.2mm, black] node (r11) at (18, 2) [circ];
\draw [line width=0.2mm, black] node (r14) at (8, 4) [amp];
\draw [line width=0.2mm, black] node (r17) at (8, 0) [amp];
\draw [line width=0.2mm, black] node (r20) at (12, 4) [amp];
\draw [line width=0.2mm, black] node (r23) at (12, 0) [amp];
\draw [line width=0.2mm, black] node (r26) at (16, 4) [amp];
\draw [line width=0.2mm, black] node (r29) at (16, 0) [amp];
\node[] at (2,4.75) {LIDG};
\node[] at (12,4.75) {RIDG};
\node[] at (2,2) {$\vec{\Psi}^T \, \vec{W}^{n+1/2}_i$};
\node[] at (8,2) {$\vec{\Psi}^T \, \widehat{\vec{W}}^{n+1/2}_{i-1}$};
\node[] at (12,2) {$\vec{\Psi}^T \, \vec{W}^{n+1/2}_i$};
\node[] at (16,2) {$\vec{\Psi}^T \, \widehat{\vec{W}}^{n+1/2}_{i+1}$};
\draw [line width=0.2mm, black] node (l1) at (2, -1.4) [amp,label=right:{upwind-in-time}];
\draw [line width=0.2mm, black] node (l2) at (7, -1.4) [circ,label=right:{interior flux}];
\draw [line width=0.2mm, black] node (l3) at (11.5, -1.4) [dia,label=right:{proper upwind flux}];
\end{tikzpicture}
\caption{Shown are the domains of dependence for the LIDG (left) and RIDG (right) prediction steps for spacetime element ${\mathcal S}^{n+1/2}_i$ in one spatial dimension.  Note that on the $t=t^n$ and $t=t^{n+1}=t^n + \Delta t$ faces, the ``proper upwind flux'' values are always on the ``past'' side of the face -- we refer to these as the ``upwind-in-time'' values. The LIDG prediction step is purely local in space -- there is no spatial communication with neighboring cells. The RIDG prediction step computes the proper upwind flux on the  $x=x_i \pm \Delta x/2$ faces, but nowhere else. In the RIDG prediction step the states to the immediate left and right of element $i$ are only temporary variables and will be discarded once the predicted solution on element $i$ has been computed -- to make note of this we place hats over the temporary variables.\label{fig:RIDG_1D}}
\end{figure}
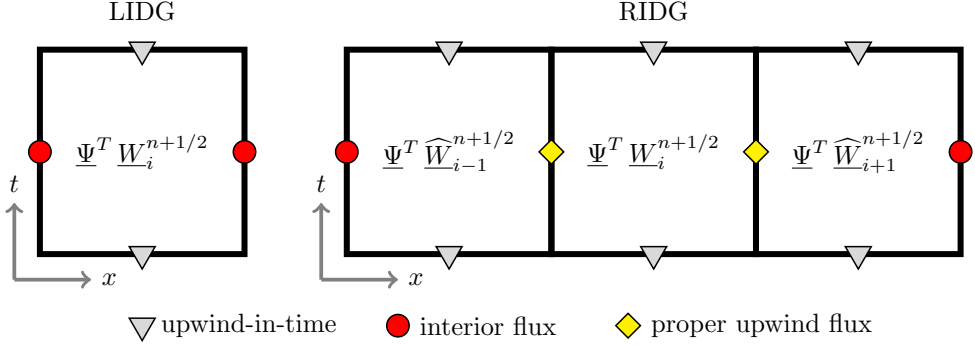

\subsection{RIDG prediction step}
\label{sec:RIDG1D_predict}
As we have argued in \cref{sec:intro} and with \cref{table:CFL_gap}, the locally-implicit predictor
described above in \cref{sec:LIDG1D_predict} will result in a scheme that has a maximum linearly stable CFL number
that is both small and becomes progressively smaller with increasing order of accuracy. We introduce a modified prediction step here to remedy these shortcomings.

The starting point is the same as in \cref{sec:LIDG1D_predict}: ansatz \cref{eqn:old_ansatz}  and 
\cref{eqn:pred_ansatz}, but now with the full spacetime ${\mathcal Q}({\mdeg},\mdim+1)$ basis in the prediction step: 
\begin{equation}
\mcorr=\mdeg+1 \qquad \text{and} \qquad \mpred=(\mdeg+1)^2.
\end{equation}
We integrate the advection equation in spacetime to get \cref{eqn:integrate_1d_adv},
but this time we integrate-by-parts in both space and time, which yields:
\begin{equation}
\label{eqn:ridg1d_pred}
\begin{split}
 \iint & \vec{\Psi} \, \left( \vec{\Psi}_{,\tau} + \nu \vec{\Psi}_{,\xi} 
 \right)^T \, \vec{W}^{n+1/2}_i \, d\tau d\xi 
+   \int  \vec{\Psi}_{|_{\tau=-1}}  \left[ \vec{\Psi}_{|_{\tau=-1}}^T \vec{W}^{n+1/2}_i
 - \vec{\Phi}^T \vec{Q}^{n}_i \right] \, d\xi \, \, - \\
 \int & \left\{ \vec{\Psi}_{|_{\xi=1}}  \left[ \nu \vec{\Psi}_{|_{\xi=1}}^T \vec{W}^{n+1/2}_i
 - {\mathcal F}^{\star}_{i+\half} \right] -  \vec{\Psi}_{|_{\xi=-1}}  \left[ \nu \vec{\Psi}_{|_{\xi=-1}}^T \vec{W}^{n+1/2}_i
 - {\mathcal F}^{\star}_{i-\half} \right] \right\} d\tau = \vec{0},
\end{split}
\end{equation}
where all 1D integrals are over $[-1,1]$, $\vec{\Phi}$ is the Legendre basis \cref{eqn:phi_basis_1d},
$\vec{\Psi}$ is the spacetime Legendre basis \cref{eqn:psi_basis_1d}, and ${\mathcal F}^{\star}$ are some appropriately defined numerical fluxes.

The crux of the idea of the regionally-implicit DG scheme in one spatial dimension can be summarized as follows:
\begin{itemize}
\item We define a {\it region} to be the current spacetime element, ${\mathcal S}^{n+1/2}_i$, and its immediate
neighbors: ${\mathcal S}^{n+1/2}_{i-1}$ and ${\mathcal S}^{n+1/2}_{i+1}$. This is illustrated in  \cref{fig:RIDG_1D}.
\item For ${\mathcal S}^{n+1/2}_i$, we use the correct upwind fluxes to define the numerical fluxes, ${\mathcal F}^{\star}$,
 on its faces.
\item For the immediate neighbors, ${\mathcal S}^{n+1/2}_{i-1}$ and ${\mathcal S}^{n+1/2}_{i+1}$, we again use the correct upwind fluxes on the faces that are shared with ${\mathcal S}^{n+1/2}_i$, but on the outer faces we use
one-sided interior fluxes. See \cref{fig:RIDG_1D}.
\item We use the ${\mathcal Q}({\mdeg},\mdim+1)$ 
spacetime basis in the prediction step (i.e., the full tensor product spacetime basis).
Numerical experimentation showed us that using the ${\mathcal Q}({\mdeg},\mdim+1)$ basis for the prediction step, rather than the
${\mathcal P}({\mdeg},\mdim+1)$ basis, produces significantly more accurate results; in the case of linear equations, this creates little additional computational expense since all the relevant matrices can be precomputed.
\end{itemize}
The result of this is a collection of three elements with solutions that are coupled to each other, but that are completely decoupled from all remaining elements. This RIDG setup is depicted in \cref{fig:RIDG_1D}, where we also show the LIDG setup as a point of comparison.

The precise form of the fluxes for the RIDG prediction step on spacetime element ${\mathcal S}^{n+1/2}_i$
can be written as follows:
\begin{gather}
\label{eqn:ridg1d_flux1}
{\mathcal F}^{\star}_{i-{3}/{2}} = \nu \vec{\Psi}_{|_{\xi=-1}}^T \vec{W}^{n+1/2}_{i-1}, \quad
{\mathcal F}^{\star}_{i-1/2} = \nu^+ \vec{\Psi}_{|_{\xi=1}}^T \vec{W}^{n+1/2}_{i-1} + \nu^- \vec{\Psi}_{|_{\xi=-1}}^T \vec{W}^{n+1/2}_{i}, \\
\label{eqn:ridg1d_flux2}
{\mathcal F}^{\star}_{i+1/2} = \nu^+ \vec{\Psi}_{|_{\xi=1}}^T \vec{W}^{n+1/2}_{i} + \nu^- \vec{\Psi}_{|_{\xi=-1}}^T \vec{W}^{n+1/2}_{i+1}, \quad
{\mathcal F}^{\star}_{i+{3}/{2}} = \nu \vec{\Psi}_{|_{\xi=1}}^T \vec{W}^{n+1/2}_{i+1}.
\end{gather}
Combining \cref{eqn:ridg1d_pred} with numerical fluxes \cref{eqn:ridg1d_flux1} and \cref{eqn:ridg1d_flux2}, yields the following
block $3 \times 3$ system:
\renewcommand{\arraystretch}{1.5}
\begin{equation}
\label{eqn:ridg1d_system}
\left[
\begin{array}{c;{2pt/2pt}c;{2pt/2pt}c}
\mat{L^0} + \mat{L^-} & \mat{X^-} & \\ \hdashline[2pt/2pt]
\mat{X^+} & \mat{L^0} + \mat{L^-} + \mat{L^+} & \mat{X^-} \\ \hdashline[2pt/2pt]
 &  \mat{X^+} & \mat{L^0} + \mat{L^+}
 \end{array}
 \right]
 \left[
\begin{array}{c}
 \widehat{\vec{W}}^{n+1/2}_{i-1} \\ \hdashline[2pt/2pt] {\vec{W}}^{n+1/2}_{i} \\ \hdashline[2pt/2pt] \widehat{\vec{W}}^{n+1/2}_{i+1}
 \end{array}
 \right] = 
 \left[
 \begin{array}{c}
 \mat{T} \, \vec{Q}^n_{i-1} \\ \hdashline[2pt/2pt] \mat{T} \, \vec{Q}^n_{i} \\ \hdashline[2pt/2pt] \mat{T} \, \vec{Q}^n_{i+1}
 \end{array}
 \right],
\end{equation}
\renewcommand{\arraystretch}{1}
where $\mat{L^0}$ is given by \cref{eqn:predictedA}, $\mat{T}$ is given by \cref{eqn:predictedB}, and
\begin{gather}
\mat{L^+} =     \frac{\nu^+}{4} \int_{-1}^{1} \vec{\Psi}_{|_{\xi=-1}} \, \vec{\Psi}_{|_{\xi=-1}}^T \, d\tau,
\quad \mat{L^-} =  -\frac{\nu^-}{4} \int_{-1}^{1} \vec{\Psi}_{|_{\xi=1}} \, \vec{\Psi}_{|_{\xi=1}}^T \, d\tau, \\
\mat{X^+} = -\frac{\nu^+}{4} \int_{-1}^{1} \vec{\Psi}_{|_{\xi=-1}} \, \vec{\Psi}_{|_{\xi=1}}^T \, d\tau, \quad
\mat{X^-} =  \frac{\nu^-}{4} \int_{-1}^{1} \vec{\Psi}_{|_{\xi=1}} \, \vec{\Psi}_{|_{\xi=-1}}^T \, d\tau,
\end{gather}
where $\mat{L^-}, \mat{L^+}, \mat{X^-}, \mat{X^+} \in \reals^{\mpred \times \mpred}$,
$\nu^+ = \max\left(\nu, 0 \right)$, and $\nu^- = \min\left(\nu, 0 \right)$.
Note that the states to the immediate left and right of the current spacetime element ${\mathcal S}^{n+1/2}_{i}$ are only temporary variables
and will be discarded once the predicted solution in element $i$ has been computed -- to make note of this we place hats over the temporary variables. This also means that we have to solve a block $3 \times 3$ system of the form \cref{eqn:ridg1d_system} on every single element ${\mathcal S}^{n+1/2}_{i}$.

\subsection{Correction step for both LIDG and RIDG}
\label{sec:RIDG1D_correct}
In order to go from the predictor to the corrector step, we multiply \cref{eqn:adv1d_nondim} by
$\vec{\Phi} \in \reals^{\mcorr}$ and integrate in spacetime:
\begin{gather}
\vec{Q}^{n+1}_i = \vec{Q}^{n}_i +  \frac{\nu}{2} \int_{-1}^{1} \int_{-1}^{1} \vec{\Phi}_{,\xi} \, q 
\, d\tau \, d\xi - \half \int_{-1}^{1} \left[ \vec{\Phi}_{|_{\xi=1}} {\mathcal F}_{i+1/2}
-  \vec{\Phi}_{|_{\xi=-1}} {\mathcal F}_{i-1/2}  \right] \, d\tau,
\end{gather}
where ${\mathcal F}$ is the numerical flux.
Next we replace $q$ by the predicted solution from either \cref{sec:LIDG1D_predict} 
or \cref{sec:RIDG1D_predict}, and use the upwind flux:
\begin{equation}
{\mathcal F}_{i-1/2} = \nu^+ \, \vec{\Psi}_{|_{\xi=1}}^T \, \vec{W}^{n+1/2}_{i-1}
+ \nu^- \, \vec{\Psi}_{|_{\xi=-1}}^T \, \vec{W}^{n+1/2}_{i},
\end{equation}
which results in
\begin{gather}
\label{eqn:correct_1d}
\vec{Q}^{n+1}_i = \vec{Q}^{n}_i + \mat{C^-} \, \vec{W}^{n+1/2}_{i-1} +  \mat{C^0} \, \vec{W}^{n+1/2}_i
+ \mat{C^+} \, \vec{W}^{n+1/2}_{i+1}, \\
\label{eqn:correct_1d_1}
\mat{C^0} = \frac{\nu}{2} \int_{-1}^{1} \int_{-1}^{1} \vec{\Phi}_{,\xi} \, \vec{\Psi}^T \, d\tau \, d\xi
- \frac{1}{2} \int_{-1}^{1} \left[ \nu^+ \vec{\Phi}_{|_{\xi=1}} \vec{\Psi}_{|_{\xi=1}}^T - 
\nu^- \vec{\Phi}_{|_{\xi=-1}} \vec{\Psi}_{|_{\xi=-1}}^T \right]  \, d\tau, \\
\label{eqn:correct_1d_2}
\mat{C^-} = \frac{\nu^+}{2} \int_{-1}^{1} \vec{\Phi}_{|_{\xi=-1}} \, \vec{\Psi}_{|_{\xi=1}}^T \, d\tau, \quad
\mat{C^+} = -\frac{\nu^-}{2} \int_{-1}^{1} \vec{\Phi}_{|_{\xi=1}} \, \vec{\Psi}_{|_{\xi=-1}}^T \, d\tau,
\end{gather}
where $\mat{C^0}, \mat{C^-}, \mat{C^+} \in \reals^{\mcorr \times \mpred}$.

\begin{figure}[!t]
\centering
(a) \includegraphics[height=45mm]{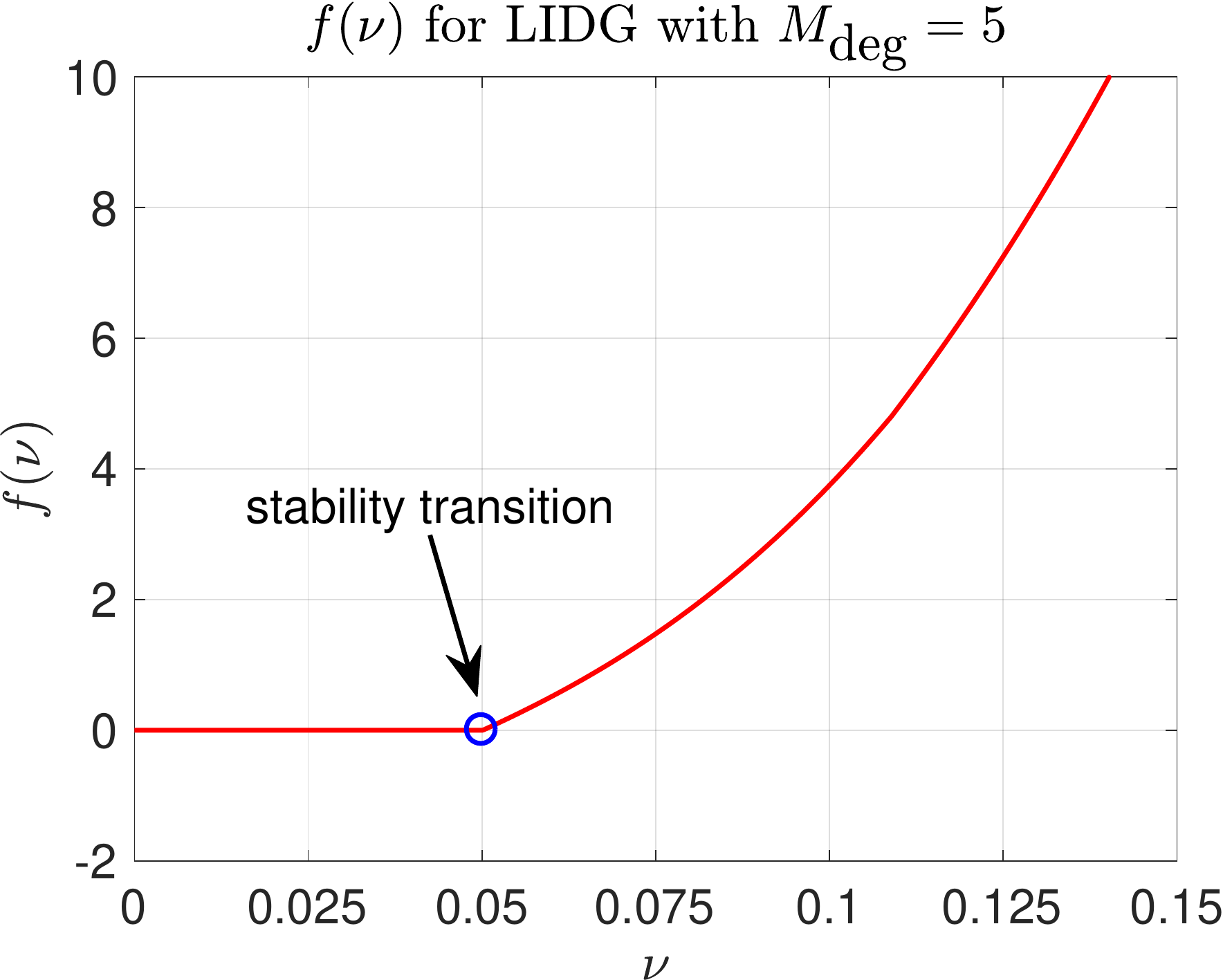} \quad
(b) \includegraphics[height=45mm]{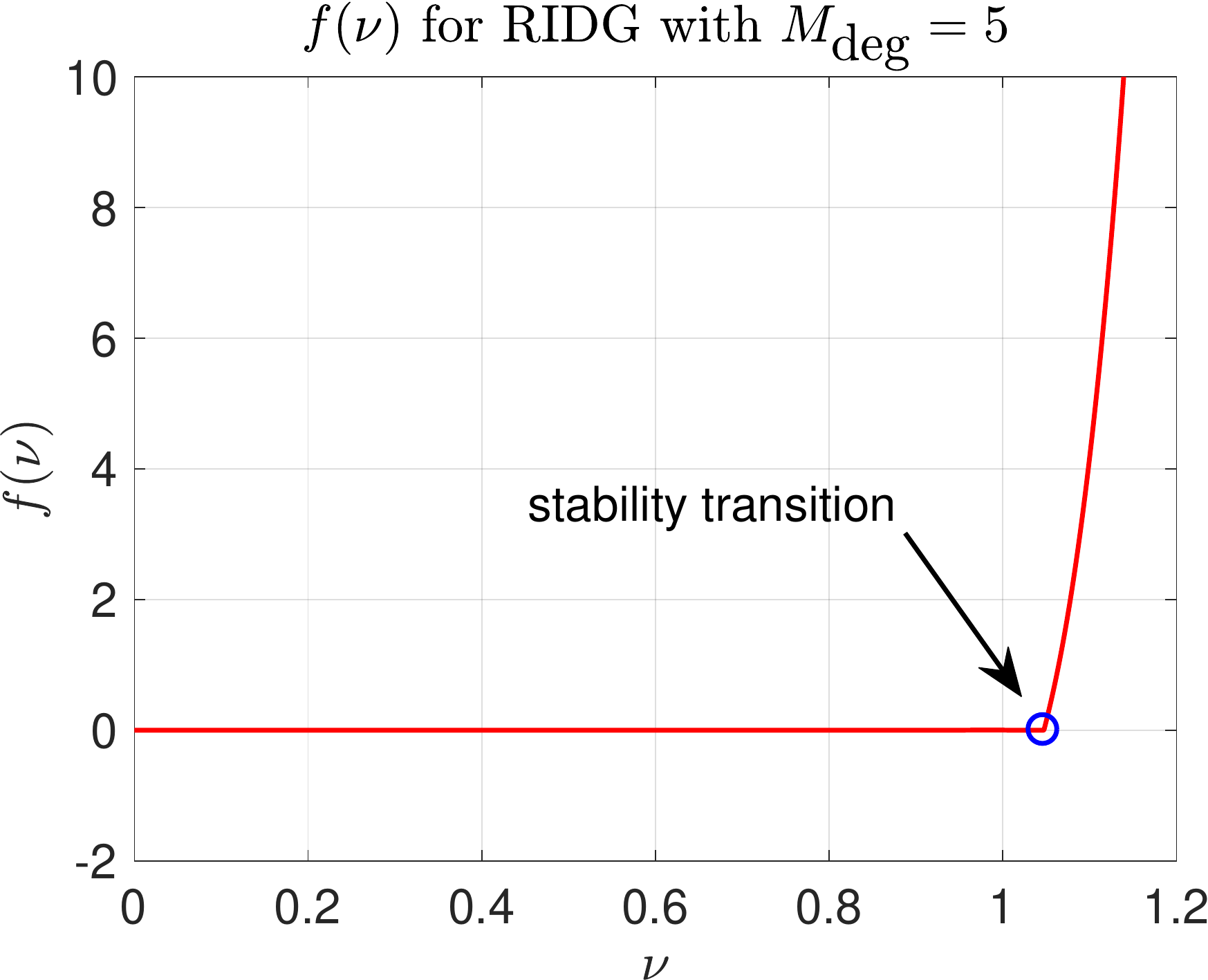}
\caption{Plot of the stability function $f(\nu)$ defined by \cref{eqn:fstab_func} for $M_{\text{deg}}=5$ for the (a) locally-implicit DG (LIDG) and (b) regionally-implicit DG (RIDG) schemes. There is a clear dichotomy between the linearly stable region and the unstable region.\label{fig:fnu}}
\end{figure}

\subsection{Von Neumann stability analysis for both LIDG and RIDG}
\label{sec:RIDG1D_stability}
The Lax-Wendroff DG scheme (aka LIDG) with prediction step detailed in \cref{sec:LIDG1D_predict} and correction step
given by \cref{eqn:correct_1d}, \cref{eqn:correct_1d_1}, and \cref{eqn:correct_1d_2}, uses a stencil involving three elements:
$\Tm_{i-1}$, $\Tm_{i}$, and $\Tm_{i+1}$; the RIDG scheme detailed in \cref{sec:RIDG1D_predict} with the same
correction step as LIDG uses a stencil involving five elements:
$\Tm_{i-2}$, $\Tm_{i-1}$, $\Tm_{i}$, $\Tm_{i+1}$, and $\Tm_{i+2}$.
This resulting CFL conditions for the LIDG and RIDG schemes are
\begin{equation}
\label{eqn:cfl_condition_1d}
| \nu | = \frac{|u| \Delta t}{\Delta x} \le 1 \qquad \text{and} \qquad
| \nu | = \frac{|u| \Delta t}{\Delta x} \le 2,
\end{equation}
repectively.
In reality, the CFL number, $|\nu|$, for which linear stability is achieved is smaller than what this CFL argument provides; we investigate this in more detail here.

In order to study linear stability,
we employ the technique of von Neumann stability analysis (e.g., see Chapter 10.5 of LeVeque \cite{book:Leveque2007}). 
In particular, we assume the following Fourier ansatz:
\begin{equation}
\label{eqn:neumann_ansatz_1d}
\vec{Q}_i^{n+1} = \vec{\widetilde{Q}}^{\, n+1} \, e^{I \omega  i} \qquad \text{and} \qquad
\vec{Q}_i^{n} = \vec{\widetilde{Q}}^{\, n} \, e^{I \omega i},
\end{equation}
where $I = \sqrt{-1}$ and $0 \le \omega \le 2\pi$ is the wave number. 
After using this ansatz, the next step is to write the resulting update in the form:
\begin{equation}
\vec{\widetilde{Q}}^{\, n+1} = \mat{\mathcal M}(\nu, \omega) \, \vec{\widetilde{Q}}^{\, n},
\end{equation}
for some matrix $\mat{\mathcal M} \in \reals^{\mcorr \times \mcorr}$.
If we apply ansatz \cref{eqn:neumann_ansatz_1d} to LIDG and RIDG, assuming w.l.o.g. that $\nu \ge 0$, we
obtain the following:
\begin{align}
\mat{{\mathcal M}_{\text{LIDG}}}(\nu,\omega) &= \left( \mat{\mathbb I} + \mat{C^0} \left(\mat{L^0}\right)^{-1} \mat{T} \right) + 
e^{-I \omega} \, \mat{C^-} \left( \mat{L^0} \right)^{-1} \mat{T}, \\
\begin{split}
\mat{{\mathcal M}_{\text{RIDG}}}(\nu,\omega) &= \left( \mat{\mathbb I} + \mat{C^0} \left( \mat{L^0} + \mat{L^+} \right)^{-1} \mat{T} \right) \\
&+ e^{-I \omega}  \left( \mat{C^-} \left( \mat{L^0} + \mat{L^+} \right)^{-1} \mat{T} - \mat{C^0} \left( \mat{L^0} + \mat{L^+} \right)^{-1}\mat{X^+} \left( \mat{L^0} \right)^{-1} \mat{T} \right) \\
&- e^{-2I\omega} \, \mat{C^{-}} \left( \mat{L^0} + \mat{L^+} \right)^{-1} \mat{X^+} \left( \mat{L^0} \right)^{-1} \mat{T},
\end{split}
\end{align}
where $\mat{\mathbb I} \in \reals^{\mcorr \times \mcorr}$ is the identity matrix.

The final step in the stability analyis is to study the spectral properties $\mat{\mathcal M}$ as a function of the CFL number $\nu$. In particular, to find the largest $\nu$ for which LIDG or RIDG are linearly stable, we define the following function:
\begin{equation}
\label{eqn:fstab_func}
f(\nu) := \max_{0 \le \omega \le 2\pi} \rho \left( \mat{\mathcal M}(\nu, \omega) \right) - 1,
\end{equation}
where $\rho\left( \mat{\mathcal M} \right)$ is the spectral radius of $\mat{\mathcal M}$. For both LIDG and RIDG,
the function $f(\nu)$ satisfies $f(0) = 0$, and there exists a finite range of $\nu$ for which
$f(\nu) \approx 0$, and there exists a value of $\nu$ for which $f(\nu)$ transitions from being approximately zero 
to rapidly increasing with increasing $\nu$. We illustrate this point in \cref{fig:fnu} for both LIDG and RIDG for the case
$\mdeg = 5$ ($\mcorr=6$, LIDG: $\mpred=21$, RIDG: $\mpred=36$); 
for each scheme we note approximately where the linear stability transition occurs.
In order to numerically estimate the location of the linear stability transition we look for the value of $\nu$ that
satisfies $f(\nu) = \varepsilon$. We do this via a simple bisection method where we set $\varepsilon = 0.0005$ and
we replace the true
maximization in \cref{eqn:fstab_func} over the maximization of 2001 uniformly spaced wave numbers over $0 \le \omega \le 2\pi$.

The result of this bisection procedure for both LIDG and RIDG for $\mdeg = 0,1,2,3,4,5$, is summarized in \cref{table:stab1d}.
In all cases we have also run the full numerical method  at various grid resolution to verify that the
simulations are indeed stable at the various CFL numbers shown in  \cref{table:stab1d}.
There are three key take-aways from  \cref{table:stab1d}:
\begin{description}
\item [\quad \textbullet] Both methods give stability regions smaller than their CFL conditions \cref{eqn:cfl_condition_1d};
\item [\quad \textbullet] The LIDG CFL number degrades roughly as the inverse of the method order; 
\item [\quad \textbullet] The RIDG CFL number is roughly one, independent of the method order. 
\end{description}

\begin{table}[!t]
\begin{center}
\begin{tabular}{|c||c|c|c|c|c|c|}
\hline 
{\bf 1D}  & $\mdeg=0$ & $\mdeg=1$ & $\mdeg=2$ & $\mdeg=3$ & $\mdeg=4$ & $\mdeg=5$ \\
\hline\hline
{\bf LIDG} & 1.000 & 0.333 & 0.171 & 0.104 & 0.070 & 0.050  \\
\hline
{\bf RIDG} & 1.000 & 1.168 & 1.135 & 1.097 & 1.066 & 1.047  \\
\hline
\end{tabular}
\caption{Numerically estimated maximum CFL numbers for the LIDG (aka Lax-Wendroff DG) and RIDG schemes in 1D.
\label{table:stab1d}}
\end{center}
\end{table}

\section{Generalization to higher dimensions}
\label{sec:higher-dimensions}
We present in this section the generalization of the proposed regionally-implicit discontinuous Galerkin (RIDG) method 
to the case of the two and three-dimensional versions of the advection equation. The key innovation beyond what was developed in \cref{sec:one-dimension} for the one-dimensional case is the inclusion of {\it transverse} cells in the prediction step.
With these inclusions, the prediction gives enhanced stability for waves propagating at all angles to the element faces.

\subsection{RIDG method in 2D}
\label{sec:ridg_in_2d}
We consider here the two-dimensional advection equation for $(t,x,y) \in [0,T] \times \Omega$ with
appropriate boundary conditions:
\begin{equation}
\label{eqn:adv2d}
q_{,t} + u_x q_{,x} + u_y q_{,y} = 0.
\end{equation}
We define a uniform Cartesian mesh with grid spacings $\Delta x$ and $\Delta y$ in each coordinate
direction. On each spacetime element:
\begin{equation}
\label{eqn:spacetime_elem_2d}
{\mathcal S}^{n+1/2}_{ij} = \left[t^{n}, t^n + \Delta t \right] \times \left[ x_i - {\Delta x}/{2}, 
x_i + {\Delta x}/{2} \right] \times \left[ y_j - {\Delta y}/{2}, 
y_j + {\Delta y}/{2} \right],
\end{equation}
we define the local coordinates, $[\tau, \xi, \eta] \in [-1,1]^3$, such that
\begin{equation}
t = t^{n+1/2} + \tau \left( {\Delta t}/{2} \right), \quad
x = x_i + \xi \left( {\Delta x}/{2} \right), \quad \text{and} \quad 
y = y_j + \eta \left( {\Delta y}/{2} \right).
\end{equation}
In these local coordinates, the advection equation is given by
\begin{equation}
q_{,\tau} + \nu_x q_{,\xi} + \nu_y q_{,\eta} = 0, \quad \nu_x = \frac{u_x \Delta t}{\Delta x}, \quad 
\nu_y = \frac{u_y \Delta t}{\Delta y},
\end{equation}
where $|\nu_x|$ and $|\nu_y|$ are the CFL numbers in each coordinate direction and the multidimensional CFL number is
\begin{equation}
\label{eqn:twod_cfl}
|\nu| := \max\left\{ |\nu_x|, |\nu_y| \right\}.
\end{equation}

At the old time, $t=t^n$, we are given the following approximate solution on each space element,
$\Tm_{ij} = \left[x_i - \Delta x/2, \, x_i + \Delta x/2\right] \times \left[y_j - \Delta y/2, \, y_j + \Delta y/2\right]$:
\begin{equation}
\label{eqn:old_ansatz_2d}
q(t^{n},x,y) \Bigl|_{{\mathcal T}_{ij}} \approx q^{n}_{ij} := 
\vec{\Phi}^T \vec{Q}^{n}_{ij},
\end{equation}
where $\vec{Q}^{n}_{ij} \in \reals^{\mcorr}$, $\vec{\Phi} \in \reals^{\mcorr}$, $\mcorr:=(\mdeg+1)(\mdeg+2)/2$, and
\begin{equation}
\label{eqn:phi_basis_2d}
\vec{\Phi} = \left( 1, \, \sqrt{3} \xi, \, \sqrt{3} \eta, \, \cdots \right),
\quad \text{s.t.} \quad \frac{1}{4} \int_{-1}^{1}\int_{-1}^{1} \vec{\Phi} \, \vec{\Phi}^T \, d\xi \, d\eta = \mat{\mathbb I} \in \reals^{\mcorr\times\mcorr},
\end{equation}
are the orthonormal space Legendre polynomials.
In order to compute a predicted solution on each spacetime element we make the following ansatz:
\begin{equation}
\label{eqn:pred_ansatz_2d}
q(t,x) \Bigl|_{{\mathcal S}^{n+1/2}_{ij}} \approx w^{n+1/2}_{ij} := 
\vec{\Psi}^T \vec{W}^{n+1/2}_{ij},
\end{equation}
where $\vec{W}^{n+1/2}_i \in \reals^{\mpred}$, $\vec{\Psi} \in \reals^{\mpred}$, $\mpred := (\mdeg+1)^3$, and
\begin{equation}
\label{eqn:psi_basis_2d}
\vec{\Psi} = \left( 1, \, \sqrt{3} \tau, \, \sqrt{3} \xi, \, \sqrt{3} \eta, \, \cdots \right), \quad
\text{s.t.} \quad \frac{1}{8} \iiint\vec{\Psi} \, \vec{\Psi}^T \, d\tau \, d\xi \, d\eta = \mat{\mathbb I} \in \reals^{\mpred\times\mpred},
\end{equation}
are the spacetime Legendre basis functions, where all 1D integrals are over $[-1,1]$.
\begin{remark}
Just as  in the one-dimensional case outlined in \cref{sec:one-dimension}, for the RIDG scheme we make use of the
${\mathcal P}({\mdeg},\mdim)$ spatial basis for the correction step (i.e., $\mcorr = (\mdeg+1)(\mdeg+2)/2$), and the
${\mathcal Q}({\mdeg},\mdim+1)$ spacetime basis for the prediction step (i.e., $\mpred = (\mdeg+1)^3$).
\end{remark}

\begin{figure}[t]
\centering
\begin{tikzpicture}[scale=0.6, 
circ/.style={scale=1.5, shape=circle, inner sep=2pt, draw, fill=red,node contents=},
dia/.style={scale=1.2, shape=diamond, inner sep=2pt, draw, fill=yellow,node contents=},
amp/.style = {scale=0.75, regular polygon, regular polygon sides=3, shape border rotate=-180,
              draw, fill=yellow, inner sep=2pt, node contents=}]

\node[] at (1,12.75) {LIDG};
\node[] at (12,12.75) {RIDG};

\draw [line width=0.8mm] (-1,0) -- (3,0) -- (3,4) -- (-1,4) -- cycle;

\draw [line width=0.2mm, black] node (r2) at (-1, 2) [circ];
\draw [line width=0.2mm, black] node (r2) at ( 3, 2) [circ];
\draw [line width=0.2mm, black] node (r2) at ( 1, 4) [circ];
\draw [line width=0.2mm, black] node (r2) at ( 1, 0) [circ];

\node[] at (1,2) {$\vec{\Psi}^T \, \vec{W}^{n+1/2}_{ij}$};

\draw[->,line width=0.5mm, gray] (5.25-7,-0.75)--(5.25-7,1) node[above,black]{$y$};
\draw[->,line width=0.5mm, gray] (5.25-7,-0.75)--(7-7,-0.75) node[right,black]{$x$};

\draw[->,line width=0.5mm, gray] (5.25,-0.75)--(5.25,1) node[above,black]{$y$};
\draw[->,line width=0.5mm, gray] (5.25,-0.75)--(7,-0.75) node[right,black]{$x$};

\draw [line width=0.8mm] (6,0) -- (10,0) -- (10,4) -- (6,4) -- cycle;
\draw [line width=0.8mm] (10,0) -- (14,0) -- (14,4) -- (10,4) -- cycle;
\draw [line width=0.8mm] (14,0) -- (18,0) -- (18,4) -- (14,4) -- cycle;

\draw [line width=0.8mm] (6,0+4) -- (10,0+4) -- (10,4+4) -- (6,4+4) -- cycle;
\draw [line width=0.8mm] (10,0+4) -- (14,0+4) -- (14,4+4) -- (10,4+4) -- cycle;
\draw [line width=0.8mm] (14,0+4) -- (18,0+4) -- (18,4+4) -- (14,4+4) -- cycle;

\draw [line width=0.8mm] (6,0+8) -- (10,0+8) -- (10,4+8) -- (6,4+8) -- cycle;
\draw [line width=0.8mm] (10,0+8) -- (14,0+8) -- (14,4+8) -- (10,4+8) -- cycle;
\draw [line width=0.8mm] (14,0+8) -- (18,0+8) -- (18,4+8) -- (14,4+8) -- cycle;

\draw [line width=0.2mm, black] node (r2) at (6, 2) [circ];
\draw [line width=0.2mm, black] node (r5) at (10, 2) [dia];
\draw [line width=0.2mm, black] node (r8) at (14, 2) [dia];
\draw [line width=0.2mm, black] node (r11) at (18, 2) [circ];
\draw [line width=0.2mm, black] node (r14) at (8, 4) [dia];
\draw [line width=0.2mm, black] node (r17) at (8, 0) [circ];
\draw [line width=0.2mm, black] node (r20) at (12, 4) [dia];
\draw [line width=0.2mm, black] node (r23) at (12, 0) [circ];
\draw [line width=0.2mm, black] node (r26) at (16, 4) [dia];
\draw [line width=0.2mm, black] node (r29) at (16, 0) [circ];
\draw [line width=0.2mm, black] node (r29) at (16, 12) [circ];
\draw [line width=0.2mm, black] node (r29) at (12, 12) [circ];
\draw [line width=0.2mm, black] node (r29) at (8, 12) [circ];
\draw [line width=0.2mm, black] node (r29) at (16, 8) [dia];
\draw [line width=0.2mm, black] node (r29) at (12, 8) [dia];
\draw [line width=0.2mm, black] node (r29) at (8, 8) [dia];
\draw [line width=0.2mm, black] node (r29) at (10, 6) [dia];
\draw [line width=0.2mm, black] node (r29) at (10, 6+4) [dia];
\draw [line width=0.2mm, black] node (r29) at (14, 6) [dia];
\draw [line width=0.2mm, black] node (r29) at (14, 6+4) [dia];
\draw [line width=0.2mm, black] node (r29) at (6, 6+4) [circ];
\draw [line width=0.2mm, black] node (r29) at (6, 6) [circ];
\draw [line width=0.2mm, black] node (r29) at (18, 6+4) [circ];
\draw [line width=0.2mm, black] node (r29) at (18, 6) [circ];

\node[] at (8,2) {$\vec{\Psi}^T \, \widehat{\vec{W}}^{n+1/2}_{i-1 \, j-1}$};
\node[] at (12,2) {$\vec{\Psi}^T \, \widehat{\vec{W}}^{n+1/2}_{i \, j-1}$};
\node[] at (16,2) {$\vec{\Psi}^T \, \widehat{\vec{W}}^{n+1/2}_{{i+1 \, j-1}}$};

\node[] at (8,2+4) {$\vec{\Psi}^T \, \widehat{\vec{W}}^{n+1/2}_{i-1 \, j}$};
\node[] at (12,2+4) {$\vec{\Psi}^T \, \vec{W}^{n+1/2}_{ij}$};
\node[] at (16,2+4) {$\vec{\Psi}^T \, \widehat{\vec{W}}^{n+1/2}_{{i+1 \, j}}$};

\node[] at (8,2+8) {$\vec{\Psi}^T \, \widehat{\vec{W}}^{n+1/2}_{i-1 \, j+1}$};
\node[] at (12,2+8) {$\vec{\Psi}^T \, \widehat{\vec{W}}^{n+1/2}_{i \, j+1}$};
\node[] at (16,2+8) {$\vec{\Psi}^T \, \widehat{\vec{W}}^{n+1/2}_{{i+1 \, j+1}}$};

\draw [line width=0.2mm, black] node (l2) at (4, -2.2) [circ,label=right:{interior flux}];
\draw [line width=0.2mm, black] node (l3) at (8.5, -2.2) [dia,label=right:{proper upwind flux}];

\end{tikzpicture}
\caption{Shown are the stencils for the LIDG (left) and RIDG (right) prediction steps in two spatial dimension. The LIDG prediction step is purely local -- there is no communication with neighboring cells. The RIDG prediction step computes the proper upwind flux on the spacetime faces shared with immediate neighbors, as well as the four corner elements. In the RIDG prediction step, all of the states, excepting only the one belonging to the middle element, are only temporary variables
and will be discarded once the predicted solution in element $ij$ has been computed -- to make note of this we place hats over the temporary variables.\label{fig:RIDG_2D}}
\end{figure}
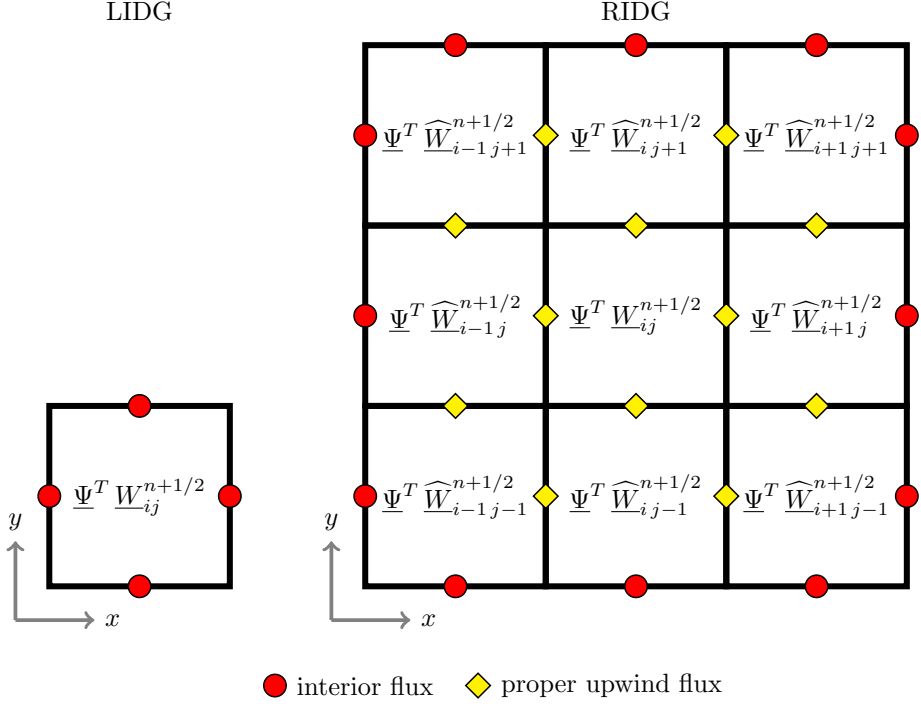

 We integrate the advection equation over a spacetime element and apply integrate-by-parts in all three
 independent variables: $\tau$, $\xi$, $\eta$, which yields: 
\begin{equation}
\label{eqn:ridg2d_pred}
\begin{split}
 & \iiint  \vec{\Psi} \,  \vec{\mathcal R}\left(\vec{\Psi}\right)^T \vec{W}^{n+1/2}_{ij}  d\vec{\mathcal S} 
+ \hspace{-1mm}  \iint \hspace{-1mm} \vec{\Psi}_{|_{\tau=-1}} \hspace{-1mm}  \left[ \vec{\Psi}_{|_{\tau=-1}}^T \vec{W}^{n+1/2}_{ij} \hspace{-1mm}
 - \vec{\Phi}^T \vec{Q}^{n}_{ij} \right] d\vec{\mathcal S}_{\tau} - \\
   \iint \hspace{-1mm} & \left\{ \vec{\Psi}_{|_{\xi=1}} \hspace{-1mm} \left[ \nu_x \vec{\Psi}_{|_{\xi=1}}^T \vec{W}^{n+1/2}_{ij}
 \hspace{-1mm} - {\mathcal F}^{\star}_{i+\half \, j} \right] \hspace{-1mm} 
   -    \vec{\Psi}_{|_{\xi=-1}} \hspace{-1mm} \left[ \nu_x \vec{\Psi}_{|_{\xi=-1}}^T \vec{W}^{n+1/2}_{ij} \hspace{-1mm}
 - {\mathcal F}^{\star}_{i-\half \, j} \right] \right\} d{\vec{\mathcal S}}_{\xi} - \\
  \iint \hspace{-1mm} & \left\{ \vec{\Psi}_{|_{\eta=1}} \hspace{-1mm} \left[ \nu_y \vec{\Psi}_{|_{\eta=1}}^T \vec{W}^{n+1/2}_{ij}
 \hspace{-1mm} - {\mathcal G}^{\star}_{i \, j+\half} \right] \hspace{-1mm} 
  -    \vec{\Psi}_{|_{\eta=-1}} \hspace{-1mm} \left[ \nu_y \vec{\Psi}_{|_{\eta=-1}}^T \vec{W}^{n+1/2}_{ij} \hspace{-1mm}
 - {\mathcal G}^{\star}_{i \, j-\half} \right] \right\} d{\vec{\mathcal S}}_{\eta} \hspace{-1mm} = \vec{0},
\end{split}
\end{equation}
where $\vec{\mathcal R}\left(\vec{\Psi}\right) = \vec{\Psi}_{,\tau} + \nu_x \vec{\Psi}_{,\xi} + \nu_y \vec{\Psi}_{,\eta}$,
$d\vec{\mathcal S} = d\tau \, d\xi \, d\eta$, $d\vec{\mathcal S}_{\tau} = d\xi \, d\eta$, 
$d\vec{\mathcal S}_{\xi} = d\tau \, d\eta$, $d\vec{\mathcal S}_{\eta} = d\tau \, d\xi$, and all 1D integrals are over $[-1,1]$.

The crux of the idea of the regionally-implicit DG scheme in two spatial dimension can be summarized as follows:
\begin{itemize}
\item We define a {\it region} to be the current spacetime element, ${\mathcal S}^{n+1/2}_{ij}$, and the eight neighbors that have a face that shares at least one point in common with one of the faces of ${\mathcal S}^{n+1/2}_{ij}$. This is illustrated in  \cref{fig:RIDG_2D}.
\item For the current spacetime element, ${\mathcal S}^{n+1/2}_{ij}$, we use the correct upwind fluxes to define
${\mathcal F}^{\star}$ and ${\mathcal G}^{\star}$.
\item For the remaining eight elements we use the correct upwind fluxes on all faces that are interior to the region and
one-sided fluxes on all faces that are on the boundary of the region. See \cref{fig:RIDG_2D}.
\item We use the ${\mathcal Q}({\mdeg},\mdim+1)$ 
spacetime basis in the prediction step (i.e., the full tensor product spacetime basis).
Numerical experimentation showed us that using the ${\mathcal Q}({\mdeg},\mdim+1)$ basis for the prediction step, rather than the
${\mathcal P}({\mdeg},\mdim+1)$ basis, produces significantly more accurate results; in the case of linear equations, this creates little additional computational expense since all the relevant matrices can be precomputed.
\end{itemize}
The result of this is a collection of nine elements with solutions that are coupled to each other, but that are completely decoupled from all remaining elements. This RIDG setup is depicted in \cref{fig:RIDG_2D}, where we also show the LIDG setup as a point of comparison.

\begin{remark}
One of the key innovations in going from the 1D RIDG scheme to its 2D counterpart is the inclusion of the {\it transverse} elements in the prediction step: ${\mathcal S}^{n+1/2}_{i-1 j-1}$, ${\mathcal S}^{n+1/2}_{i+1 j-1}$, ${\mathcal S}^{n+1/2}_{i-1 j+1}$, 
and ${\mathcal S}^{n+1/2}_{i+1 j+1}$. Without these transverse cells, the maximum allowable two-dimensional CFL number \cref{eqn:twod_cfl} remains small for any waves traveling transverse to the mesh. In the recent literature, there exist several variants of genuinely multidimensional Riemann solvers (e.g., Balsara \cite{article:Balsara2014}); by including the transverse elements, the current work can be viewed as an example of a novel type of multidimensional Riemann solver.
\end{remark}

Applying all of the above principles to \cref{eqn:ridg2d_pred} for all of the nine elements that are in the current region
yields a block  $9 \times 9$ linear system. The left-hand side of this system can be written as
\renewcommand{\arraystretch}{1.5}
\begin{equation}
\label{eqn:ridg2d_system}
\left[
\hspace{-1mm}
\begin{array}{c;{2pt/2pt}c;{2pt/2pt}c;{2pt/2pt}c;{2pt/2pt}c;{2pt/2pt}c;{2pt/2pt}c;{2pt/2pt}c;{2pt/2pt}c}
\Lone & \mat{X^-} & & \mat{Y^-} & & & & &\\ \hdashline[2pt/2pt]
\mat{X^+} & \Ltwo & \mat{X^-} & & \mat{Y^-} & & & &\\ \hdashline[2pt/2pt]
& \mat{X^+} & \Lthree & & & \mat{Y^-} & & &\\ \hdashline[2pt/2pt]
\mat{Y^+} & & & \Lfour & \mat{X^-} & & \mat{Y^-} & &\\ \hdashline[2pt/2pt]
& \mat{Y^+} & & \mat{X^+} & \Lfive & \mat{X^-} & & \mat{Y^-} &\\ \hdashline[2pt/2pt]
& & \mat{Y^+} & & \mat{X^+} & \Lsix & & & \mat{Y^-} \\ \hdashline[2pt/2pt]
& & & \mat{Y^+} & & & \Lseven & \mat{X^-} &\\ \hdashline[2pt/2pt]
& & & & \mat{Y^+} & & \mat{X^+} & \Leight & \mat{X^-}\\ \hdashline[2pt/2pt]
& & & & & \mat{Y^+} & & \mat{X^+} & \Lnine
 \end{array}
 \hspace{-1mm}
 \right] \hspace{-2mm}
 \left[ \hspace{-1mm}
\begin{array}{c}
 \widehat{\vec{W}}^{n+1/2}_{i-1  j-1} \\ \hdashline[2pt/2pt] \widehat{\vec{W}}^{n+1/2}_{i  j-1} \\ \hdashline[2pt/2pt] \widehat{\vec{W}}^{n+1/2}_{i+1  j-1} \\ \hdashline[2pt/2pt]
 \widehat{\vec{W}}^{n+1/2}_{i-1  j} \\ \hdashline[2pt/2pt] {\vec{W}}^{n+1/2}_{i  j} \\ \hdashline[2pt/2pt] \widehat{\vec{W}}^{n+1/2}_{i+1  j} \\ \hdashline[2pt/2pt]
 \widehat{\vec{W}}^{n+1/2}_{i-1  j+1} \\ \hdashline[2pt/2pt] \widehat{\vec{W}}^{n+1/2}_{i  j+1} \\ \hdashline[2pt/2pt] \widehat{\vec{W}}^{n+1/2}_{i+1  j+1}
 \end{array}
 \hspace{-1mm}
 \right],
\end{equation}
\renewcommand{\arraystretch}{1}
and the right-hand side can be written as
\begin{equation}
\begin{split}
\Bigl[ \, \,
&\mat{T} \, \vec{Q}^{n}_{i-1  j-1}, \quad 
\mat{T} \, \vec{Q}^{n}_{i  j-1}, \quad \cdots, \quad
\mat{T} \, \vec{Q}^{n}_{i  j}, \quad \cdots, \quad
\mat{T} \, \vec{Q}^{n}_{i  j+1}, \quad
\mat{T} \, \vec{Q}^{n}_{i+1  j+1} \, \Bigr]^T,
\end{split}
\end{equation}
where
\begin{gather}
\mat{X}^{\pm} = \mp \frac{\nu_x^\pm}{8} \iint \vec{\Psi}_{|_{\xi=\mp1}} \, \vec{\Psi}_{|_{\xi=\pm1}}^T \, d\vec{\mathcal S}_{\xi}, \quad
\mat{Y}^{\pm} = \mp \frac{\nu_y^\pm}{8} \iint \vec{\Psi}_{|_{\eta=\mp1}} \, \vec{\Psi}_{|_{\eta=\pm1}}^T \, d\vec{\mathcal S}_{\eta}, \\
\mat{L^{\alpha \beta \gamma \delta}} = \mat{L^{0}} + \alpha \, \mat{L^{-}_x} + \beta \, \mat{L^{+}_x} + \gamma \, \mat{L^{-}_y} + 
\delta \, \mat{L^{+}_y}, \quad \alpha, \beta, \gamma, \delta \in \left\{ 0 , 1 \right\}, \\
\mat{L^0} =  \frac{1}{8} \iiint \vec{\Psi} \, \vec{\mathcal R}\left(\vec{\Psi}\right)^T \hspace{-1mm} d\vec{\mathcal S} 
+  \frac{1}{8} \iint \vec{\Psi}_{|_{\tau=-1}} \, \vec{\Psi}_{|_{\tau=-1}}^T  d\vec{\mathcal S}_{\tau}, \quad
\mat{T} = \frac{1}{8}  \iint \vec{\Psi}_{|_{\tau=-1}} \, \vec{\Phi}^T   d\vec{\mathcal S}_{\tau}, \\
\mat{L_x^\pm} =   \pm  \frac{\nu_x^{\pm}}{8} \iint \vec{\Psi}_{|_{\xi=\mp1}} \, \vec{\Psi}_{|_{\xi=\mp1}}^T \, d\vec{\mathcal S}_{\xi},
\quad
\mat{L_y^\pm} =   \pm  \frac{\nu_y^{\pm}}{8} \iint \vec{\Psi}_{|_{\eta=\mp1}} \, \vec{\Psi}_{|_{\eta=\mp1}}^T \, d\vec{\mathcal S}_{\eta},
\end{gather}
where $\mat{T} \in \reals^{\mpred \times \mcorr}$ and $\mat{X}^{\pm}, \mat{Y}^{\pm}, \mat{L^0}, \mat{L_x^\pm}, \mat{L_y^\pm}
\in \reals^{\mpred \times \mpred}$.

The correction step can be written as
\begin{gather}
\label{eqn:correct_2d}
\vec{Q}^{n+1}_{ij} = \vec{Q}^{n}_{ij} + \mat{C_x^-} \, \vec{W}^{n+\half}_{i-1 j} +
\mat{C_y^-} \, \vec{W}^{n+\half}_{i j-1} +  \mat{C^0} \, \vec{W}^{n+\half}_{ij}
+ \mat{C_x^+} \, \vec{W}^{n+\half}_{i+1 j}
+  \mat{C_y^+} \, \vec{W}^{n+\half}_{ij+1}, \\
\label{eqn:correct_2d_1}
\begin{split}
\mat{C^0} = \frac{1}{4} \iiint \vec{\mathcal U}\left( \vec{\Phi} \right) \, \vec{\Psi}^T d\vec{\mathcal S} 
- \frac{1}{4} & \iint \left[ \nu^+_x \vec{\Phi}_{|_{\xi=1}} \vec{\Psi}_{|_{\xi=1}}^T - 
\nu^-_x \vec{\Phi}_{|_{\xi=-1}} \vec{\Psi}_{|_{\xi=-1}}^T \right]   d\vec{\mathcal S}_{\xi} \\
- \frac{1}{4} & \iint \left[ \nu^+_y \vec{\Phi}_{|_{\eta=1}} \vec{\Psi}_{|_{\eta=1}}^T - 
\nu^-_y \vec{\Phi}_{|_{\eta=-1}} \vec{\Psi}_{|_{\eta=-1}}^T \right]   d\vec{\mathcal S}_{\eta}, 
\end{split} \\
\label{eqn:correct_2d_2}
\mat{C^{\mp}_x} = \pm \frac{\nu^{\pm}_x}{4} \iint \vec{\Phi}_{|_{\xi=\mp1}} \, \vec{\Psi}_{|_{\xi=\pm1}}^T d\vec{\mathcal S}_{\xi}, \quad
\mat{C^{\mp}_y} = \pm \frac{\nu^{\pm}_y}{4} \iint \vec{\Phi}_{|_{\eta=\mp1}} \, \vec{\Psi}_{|_{\eta=\pm1}}^T d\vec{\mathcal S}_{\eta},
\end{gather}
where \, $\vec{\mathcal U}(\vec{\Phi}) = \nu_x \vec{\Phi}_{,\xi} + \nu_y \vec{\Phi}_{,\eta}$, \,
$\mat{C^0}, \mat{C^\pm_x}, \mat{C^\pm_y} \in \reals^{\mcorr \times \mpred}$.

\subsection{RIDG method in 3D}
We consider here the three-dimensional advection equation for $(t,x,y,z) \in [0,T] \times \Omega$ with
appropriate boundary conditions:
\begin{equation}
\label{eqn:adv3d}
q_{,t} + u_x q_{,x} + u_y q_{,y} + u_z q_{,z} = 0.
\end{equation}
We define a uniform Cartesian mesh with grid spacings $\Delta x$, $\Delta y$, and $\Delta z$ in each coordinate
direction. On each spacetime element:
\begin{equation}
\label{eqn:spacetime_elem_3d}
{\mathcal S}^{n+1/2}_{ijk} =  {\mathcal S}^{n+1/2}_{ij} \times \left[ z_k - {\Delta z}/{2}, 
z_k + {\Delta z}/{2} \right],
\end{equation}
where ${\mathcal S}^{n+1/2}_{ij}$ is defined by \cref{eqn:spacetime_elem_2d}, 
we define the local coordinates, $[\tau, \xi, \eta, \zeta] \in [-1,1]^4$, such that
\begin{equation}
t = t^{n+1/2} + \tau \left( {\Delta t}/{2} \right), \,
x = x_i + \xi \left( {\Delta x}/{2} \right), \,
y = y_j + \eta \left( {\Delta y}/{2} \right), \,  
z = z_k + \zeta \left( {\Delta z}/{2} \right).
\end{equation}
In these local coordinates, the advection equation is given by
\begin{equation}
q_{,\tau} + \nu_x q_{,\xi} + \nu_y q_{,\eta} + \nu_z q_{,\zeta} = 0, 
\quad \nu_x = \frac{u_x \Delta t}{\Delta x}, \quad 
\nu_y = \frac{u_y \Delta t}{\Delta y}, \quad
\nu_z = \frac{u_z \Delta t}{\Delta z},
\end{equation}
where $|\nu_x|$, $|\nu_y|$, and $|\nu_z|$ are the CFL numbers in each coordinate direction and the multidimensional CFL number is
\begin{equation}
\label{eqn:threed_cfl}
|\nu| := \max\left\{ |\nu_x|, |\nu_y|, |\nu_z| \right\}.
\end{equation}

The development of the RIDG scheme in 3D is completely analogous to the 2D RIDG scheme from \cref{sec:ridg_in_2d}. In 1D the prediction step requires a stencil of 3 elements, in 2D we need $3^2=9$ elements, and in
3D we need $3^3 = 27$ elements. For the sake of brevity we omit the details.

\begin{table}[!t]
\begin{center}
\begin{tabular}{|c||c|c|c|c|c|c|}
\hline 
 {\bf 2D} & $\mdeg=0$ & $\mdeg=1$ & $\mdeg=3$ & $\mdeg=5$ & $\mdeg=7$ & $\mdeg=9$ \\
\hline\hline
{\bf LIDG} & 0.50 & 0.23 & 0.08 & 0.04 & 0.025 & 0.01  \\
\hline
{\bf RIDG} & 1.00 & 1.00 & 0.80 & 0.75 & 0.75 & 0.75  \\
\hline \hline
{\bf 3D} & $\mdeg=0$ & $\mdeg=1$ & $\mdeg=3$ & $\mdeg=5$ & $\mdeg=7$ & $\mdeg=9$ \\
\hline\hline
{\bf LIDG} & 0.33 & 0.10 & 0.03 & 0.025 & 0.02 & 0.01  \\
\hline
{\bf RIDG} & 1.00 & 0.80  & 0.60 & 0.60 & 0.60 & 0.60  \\
\hline
\end{tabular}
\caption{Numerically estimated maximum CFL numbers, $|\nu|$, for the LIDG (aka Lax-Wendroff DG) and RIDG schemes in 2D and 3D.
\label{table:stab2d}}
\end{center}
\end{table}

\subsection{Von Neumann stability analysis for both LIDG and RIDG}
Linear stability analysis proceeds in 2D and 3D in the same manner as in 1D. We take the numerical update and make the
Fourier ansatz:
\begin{equation}
\vec{Q}_{ijk}^{n+1} = \vec{\widetilde{Q}}^{\, n+1} \, e^{I \left( \omega_x i + \omega_y j + \omega_z k \right)} \qquad \text{and} \qquad
\vec{Q}_{ijk}^{n} = \vec{\widetilde{Q}}^{\, n} \, e^{I \left( \omega_x i + \omega_y j + \omega_z k \right)},
\end{equation}
where $I = \sqrt{-1}$ and $0 \le \omega_x, \, \omega_y, \, \omega_z \le 2\pi$ are the wave numbers in each coordinate direction.
After using this ansatz, we again write the resulting update in the form:
\begin{equation}
\vec{\widetilde{Q}}^{\, n+1} = \mat{\mathcal M}\left(\nu_x, \nu_y, \nu_z, \omega_x, \omega_y, \omega_z \right) 
\, \vec{\widetilde{Q}}^{\, n},
\end{equation}
for some matrix $\mat{\mathcal M} \in \reals^{\mcorr \times \mcorr}$.
Finally, we define the function
\begin{equation}
\label{eqn:fstab_func_3d}
f(\nu_x, \nu_y, \nu_z) := 
\max_{0 \le \omega_x, \, \omega_y, \, \omega_z \le 2\pi}  \rho \left( \mat{\mathcal M}(\nu_x, \, \nu_y, \, \nu_z, \, \omega_x, \, \omega_y, \, \omega_z) \right) - 1,
\end{equation}
where $\rho\left( \mat{\mathcal M} \right)$ is the spectral radius of $\mat{\mathcal M}$.

\begin{figure}[!th]
\begin{tabular}{cc}
2D LIDG $\mdeg=1$ & 2D RIDG $\mdeg=1$ \\
(a)\includegraphics[width = 0.43\textwidth]{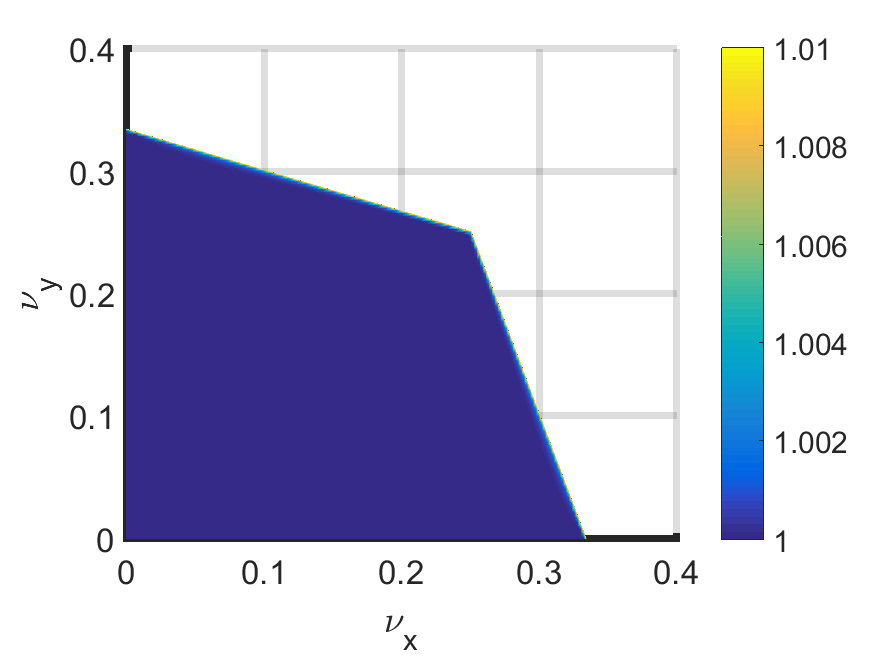} &
(b)\includegraphics[width = 0.43\textwidth]{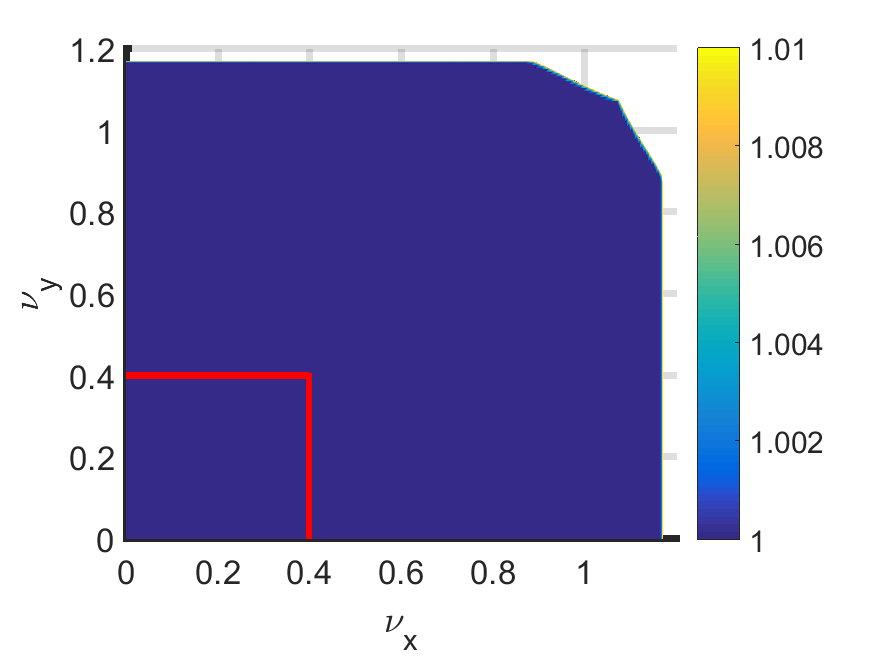} \\
2D LIDG $\mdeg=3$ & 2D RIDG $\mdeg=3$ \\
(c)\includegraphics[width = 0.43\textwidth]{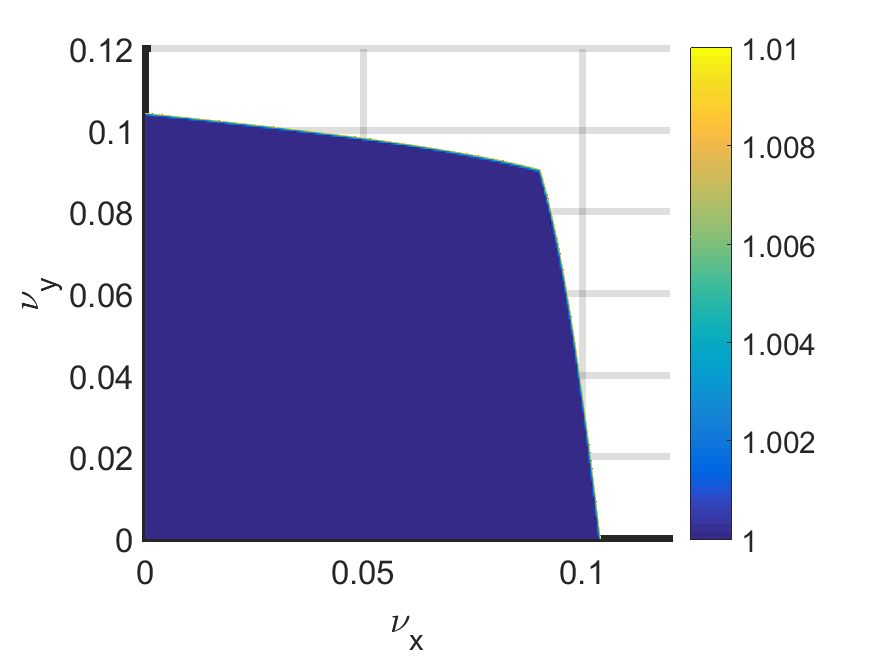} &
(d)\includegraphics[width = 0.43\textwidth]{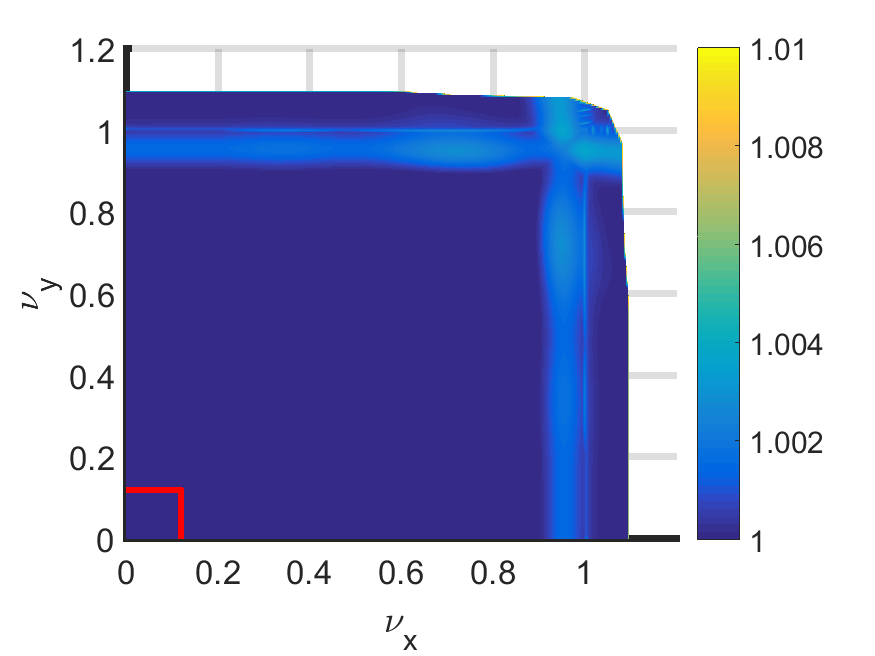} \\
2D LIDG $\mdeg=5$ & 2D RIDG $\mdeg=5$ \\
(e)\includegraphics[width = 0.43\textwidth]{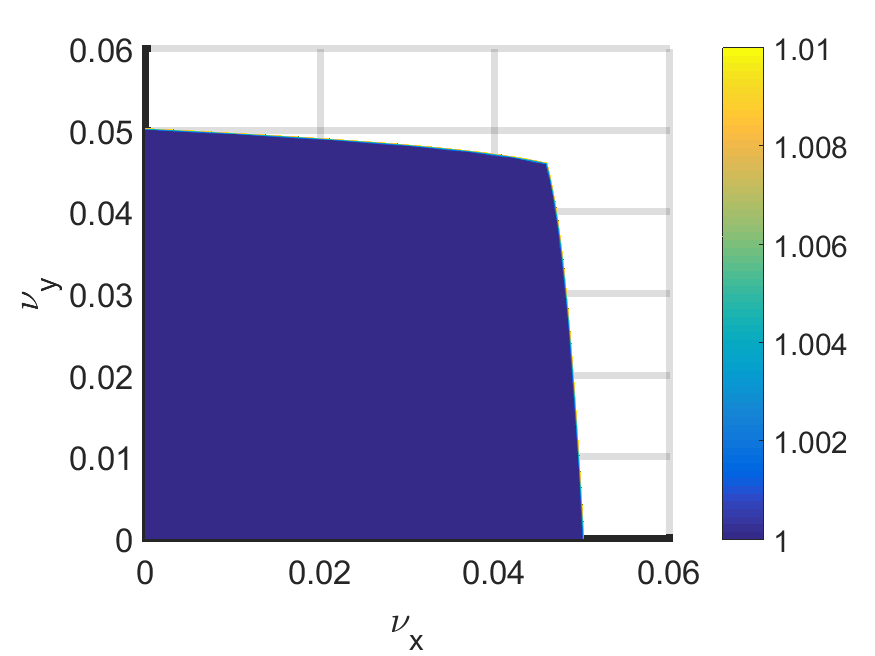} &
(f)\includegraphics[width = 0.43\textwidth]{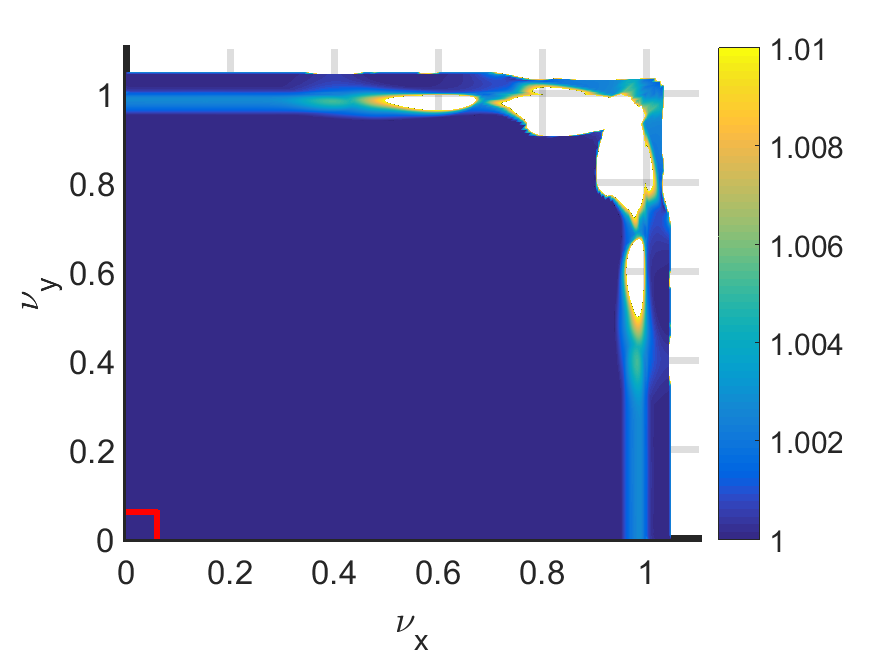}  \\
2D LIDG $\mdeg=7$ & 2D RIDG $\mdeg=7$ \\
(g)\includegraphics[width = 0.43\textwidth]{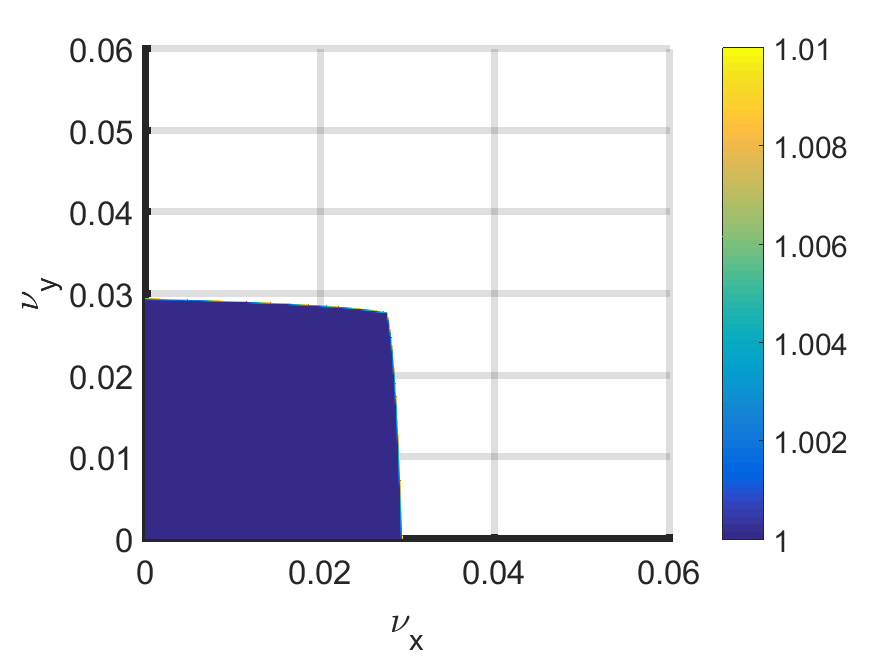} &
(h)\includegraphics[width = 0.43\textwidth]{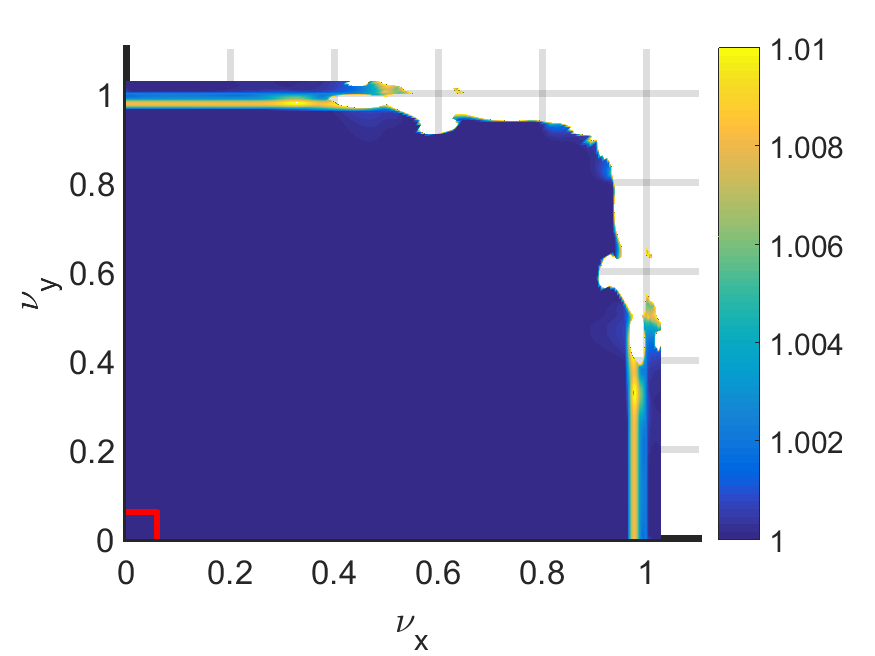}
\end{tabular}
\caption{Stability plots for the two-dimensional LIDG and RIDG methods for various polynomial orders. 
Shown is a false color plot of $f(\nu_x,\nu_y)+1$ as defined
by \cref{eqn:fstab_func_3d}. Note the different horizontal and vertical scales on the plots of the LIDG and RIDG schemes.
The red box in each RIDG plot demonstrates the plotting bounds for the LIDG stability region for the same method order. 
From these plots we can estimate the maximum CFL number $|\nu|$ as defined by \cref{eqn:twod_cfl}: 
(a) $|\nu| \lessapprox 0.23$ ($\mdeg=1$, LIDG), (b) $|\nu| \lessapprox 1.00$ ($\mdeg=1$, RIDG), 
(c) $|\nu| \lessapprox 0.08$ ($\mdeg=3$, LIDG), (d) $|\nu| \lessapprox 0.80$ ($\mdeg=3$, RIDG), 
(e) $|\nu| \lessapprox 0.04$ ($\mdeg=5$, LIDG), (f) $|\nu| \lessapprox 0.75$ ($\mdeg=5$, RIDG), 
(g) $|\nu| \lessapprox 0.025$ ($\mdeg=7$, LIDG), and (h) $|\nu| \lessapprox 0.75$ ($\mdeg=7$, RIDG).
\label{fig:stabilityr0m2}}
\end{figure}

Just as in 1D, we estimate the maximum CFL numbers of LIDG and RIDG by studying the values of \cref{eqn:fstab_func_3d}.
Our numerically obtained estimates for the maximum value of $|\nu|$ as defined by \cref{eqn:twod_cfl} and 
\cref{eqn:threed_cfl} are summarized in \cref{table:stab2d}. 
Again we see the following:
\begin{itemize}
\item LIDG: the maximum stable CFL number tends to zero as the polynomial degree is increased; and
\item RIDG: the maximum stable CFL number has a finite lower bound with increasing polynomial degree (approximately $0.75$ in 2D and $0.60$ in 3D).
\end{itemize}
To get a more detailed view of the stability function \cref{eqn:fstab_func_3d} in 2D, 
we show false color plots of $f(\nu_x,\nu_y)+1$ in \cref{fig:stabilityr0m2} for both LIDG and RIDG for various method orders.
The transverse elements that were included in the prediction step for RIDG (see \cref{sec:ridg_in_2d}) are critically
important in achieving a stability region that does not significantly degrade in going from 1D to 2D.

\section{Numerical convergence studies}
\label{sec:results}
In this section we present convergence studies in 1D, 2D, and 3D for both LIDG and RIDG, and compare the errors and runtimes for the two methods. In all cases, we compute an approximate order of accuracy using the following approximation:
\begin{equation}
\label{eqn:Mratio}
\text{error}(h) = c h^M + {\mathcal O}\left(h^{M+1} \right) \quad \Longrightarrow \quad
M \approx \frac{\log\paren{{\text{error}(h_1)}/{\text{error}(h_2)}}}{\log\paren{{h_1}/{h_2}} }.
\end{equation}


\begin{table}[!t]
\centering
\begin{tabular}{|r|r|c|c|c|c|c|c|}
\hline	
\multicolumn{8}{|c|}{{\bf 1D LIDG: ($\mdeg=3$, $\nu=0.104$)}}	\\ \hline
{\bf mesh}		&	{\bf $T_r$(s)}	&	$L^1$ {\bf error}	& \cref{eqn:Mratio}	&	$L^2$ {\bf error}	& \cref{eqn:Mratio}
	&	$L^\infty$ {\bf error}	& \cref{eqn:Mratio}	\\ \hline
$	40	$&	0.221	&$	1.83	e{	-1	}$&	--	&$	1.83	e{	-1	}$&	--	&$	1.92	e{	-1	}$&	--	\\			
$	80	$&	0.784	&$	1.08	e{	-2	}$&	4.08	&$	1.07	e{	-2	}$&	4.09	&$	1.13	e{	-2	}$&	4.09	\\	
$	160	$&	3.073	&$	6.52	e{	-4	}$&	4.05	&$	6.46	e{	-4	}$&	4.05	&$	6.66	e{	-4	}$&	4.09	\\	
$	320	$&	12.246	&$	4.01	e{	-5	}$&	4.02	&$	4.00	e{	-5	}$&	4.01	&$	4.10	e{	-5	}$&	4.02	\\			
$	640	$&	48.928	&$	2.49	e{	-6	}$&	4.01	&$	2.50	e{	-6	}$&	4.00	&$	2.79	e{	-6	}$&	3.88	\\			\hline\hline
\multicolumn{8}{|c|}{{\bf 1D RIDG: ($\mdeg=3$, $\nu=0.9$)}}	\\ \hline		
{\bf mesh}		&	{\bf $T_r$(s)}	&	$L^1$ {\bf error}	& \cref{eqn:Mratio}	&	$L^2$ {\bf error}	& \cref{eqn:Mratio}
	&	$L^\infty$ {\bf error}	& \cref{eqn:Mratio}	\\ \hline
$	40	$&	0.032	&$	8.46	e{	-2	}$&	--	&$	8.77	e{	-2	}$&	--	&$	1.02	e{	-1	}$&	--	\\			
$	80	$&	0.121	&$	3.67	e{	-3	}$&	4.53	&$	3.72	e{	-3	}$&	4.56	&$	4.68	e{	-3	}$&	4.45	\\	
$	160	$&	0.473	&$	1.51	e{	-4	}$&	4.61	&$	1.52	e{	-4	}$&	4.62	&$	1.76	e{	-4	}$&	4.73	\\	
$	320	$&	1.885	&$	7.96	e{	-6	}$&	4.24	&$	8.02	e{	-6	}$&	4.24	&$	8.95	e{	-6	}$&	4.30	\\		
$	640	$&	7.618	&$	4.75	e{	-7	}$&	4.07	&$	4.77	e{	-7	}$&	4.07	&$	5.57	e{	-7	}$&	4.01	\\			\hline\hline
\multicolumn{8}{|c|}{{\bf 1D LIDG: ($\mdeg=5$, $\nu=0.04$)}}	\\ \hline												
{\bf mesh}		&	{\bf $T_r$(s)}	&	$L^1$ {\bf error}	& \cref{eqn:Mratio}	&	$L^2$ {\bf error}	& \cref{eqn:Mratio}
	&	$L^\infty$ {\bf error}	& \cref{eqn:Mratio}	\\ \hline
$	40	$&	0.518	&$	1.11	e{	-3	}$&	--	&$	1.11	e{	-3	}$&	--	&$	1.25	e{	-3	}$&	--	\\			
$	80	$&	2.041	&$	1.74	e{	-5	}$&	6.00	&$	1.76	e{	-5	}$&	5.98	&$	1.88	e{	-5	}$&	6.06	\\		
$	160	$&	8.060	&$	2.73	e{	-7	}$&	5.99	&$	2.72	e{	-7	}$&	6.02	&$	2.86	e{	-7	}$&	6.04	\\	
$	320	$&	32.196	&$	4.24	e{	-9	}$&	6.01	&$	4.23	e{	-9	}$&	6.01	&$	4.36	e{	-9	}$&	6.03	\\	
$	640	$&	129.172	&$	6.61	e{	-11	}$&	6.00	&$	6.61	e{	-11	}$&	6.00	&$	6.78	e{	-11	}$&	6.01	\\			\hline\hline
\multicolumn{8}{|c|}{{\bf 1D RIDG: ($\mdeg=5$, $\nu=0.9$)}}	\\ \hline
{\bf mesh}		&	{\bf $T_r$ (s)}	&	$L^1$ {\bf error}	& \cref{eqn:Mratio}	&	$L^2$ {\bf error}	& \cref{eqn:Mratio}
	&	$L^\infty$ {\bf error}	& \cref{eqn:Mratio}	\\ \hline
$	40	$&	0.033	&$	1.50	e{	-4	}$&	--	&$	1.65	e{	-4	}$&	--	&$	4.64	e{	-4	}$&	--	\\			
$	80	$&	0.123	&$	2.68	e{	-6	}$&	5.81	&$	2.79	e{	-6	}$&	5.89	&$	5.19	e{	-6	}$&	6.48	\\	
$	160	$&	0.482	&$	3.91	e{	-8	}$&	6.10	&$	4.05	e{	-8	}$&	6.11	&$	4.89	e{	-8	}$&	6.73	\\		
$	320	$&	1.936	&$	5.85	e{	-10	}$&	6.06	&$	6.12	e{	-10	}$&	6.05	&$	8.37	e{	-10	}$&	5.87	\\	
$	640	$&	7.733	&$	8.94	e{	-12	}$&	6.03	&$	9.46	e{	-12	}$&	6.02	&$	1.36	e{	-11	}$&	5.94	\\			\hline
\end{tabular}
\caption{Convergence and runtime study for the 1D LIDG and RIDG  methods with $\mdeg=3$ and $\mdeg=5$
on 1D advection equation \cref{eqn:adv1d} with initial condition \cref{cc1D}.  Shown for various mesh sizes are the runtimes 
($T_r$) measured in seconds, 
the relative errors in $L^1$, $L^2$, and $L^{\infty}$, as well as the estimated convergence rates according to formula \cref{eqn:Mratio}.  We see that for any fixed number of elements, the RIDG method has a shorter runtime (i.e., computational cost) and smaller error.
\label{table:RIvLI1D}}
\end{table}%


\subsection{1D convergence tests}
We consider the 1D advection equation \cref{eqn:adv1d} with $u=1$, $\Omega=[-1,1]$, periodic BCs,
and initial condition: 
\begin{equation} 
\label{cc1D}
q(t=0,x) = \sin(16 \pi x).
\end{equation}
We run the code \cite{code:ridg-code} to $t=2$ with $\mdeg=3$ (LIDG: $\nu=0.104$, RIDG: $\nu=0.9$) and $\mdeg=5$ (LIDG: $\nu=0.04$, RIDG: $\nu=0.9$) and compare runtimes and errors;  the results are shown in \cref{table:RIvLI1D}. 
We see that both methods exhibit the expected convergence rates in 
$L^1$, $L^2$, and $L^\infty$.
For a fixed number of elements, usage of the RIDG method leads to smaller errors.   We also notice that for a fixed number of elements the experiment runtime for the RIDG method is shorter than that of the LIDG method -- this is due to the increase in the maximum linearly stable CFL number from LIDG to RIDG.


\begin{table}[!t]		
\centering
\begin{tabular}{|r|r|c|c|c|c|c|c|}
\hline	
\multicolumn{8}{|c|}{{\bf 2D LIDG: ($\mdeg=3$, $\nu=0.05$)}}	\\ \hline
{\bf mesh}		&	{\bf$T_r$ (s)}	&	$L^1$ {\bf error}	& \cref{eqn:Mratio}	&	$L^2$ {\bf error}	& \cref{eqn:Mratio}
	&	$L^\infty$ {\bf error}	& \cref{eqn:Mratio}	\\ \hline																					
$	40	^2$&	22.2	&$	8.75	e{	-1	}$&	--	&$	7.87	e{	-1	}$&	--	&$	7.93	e{	-1	}$&	--	\\			
$	80	^2$&	176.8	&$	6.37	e{	-2	}$&	3.78	&$	5.72	e{	-2	}$&	3.78	&$	6.54	e{	-2	}$&	3.6	\\			
$	160	^2$&	1415.7	&$	1.98	e{	-3	}$&	5.01	&$	1.81	e{	-3	}$&	4.98	&$	2.94	e{	-3	}$&	4.48	\\			
$	320	^2$&	11335.2	&$	7.66	e{	-5	}$&	4.69	&$	7.09	e{	-5	}$&	4.67	&$	1.61	e{	-4	}$&	4.19	\\			
\hline\hline																									
\multicolumn{8}{|c|}{{\bf 2D RIDG: ($\mdeg=3$, $\nu=0.75$)}}	\\ \hline																						
{\bf mesh}		&	{\bf$T_r $ (s)}	&	$L^1$ {\bf error}	& \cref{eqn:Mratio}	&	$L^2$ {\bf error}	& \cref{eqn:Mratio}										
	&	$L^\infty$ {\bf error}	& \cref{eqn:Mratio}	\\ \hline																					
$	40	^2$&	3.5	&$	6.29	e{	-1	}$&	--	&$	5.58	e{	-1	}$&	--	&$	5.62	e{	-1	}$&	--	\\			
$	80	^2$&	27.1	&$	2.81	e{	-2	}$&	4.49	&$	2.54	e{	-2	}$&	4.46	&$	3.45	e{	-2	}$&	4.03	\\			
$	160	^2$&	218.6	&$	1.04	e{	-3	}$&	4.75	&$	9.58	e{	-4	}$&	4.73	&$	1.76	e{	-3	}$&	4.29	\\			
$	320	^2$&	1740.0	&$	5.37	e{	-5	}$&	4.28	&$	4.95	e{	-5	}$&	4.27	&$	1.09	e{	-4	}$&	4.02	\\			
\hline\hline																									
\multicolumn{8}{|c|}{{\bf 2D LIDG: ($\mdeg=5$, $\nu=0.03$)}}	\\ \hline																								
{\bf mesh}		&	{\bf$T_r $ (s)}	&	$L^1$ {\bf error}	& \cref{eqn:Mratio}	&	$L^2$ {\bf error}	& \cref{eqn:Mratio}																
	&	$L^\infty$ {\bf error}	& \cref{eqn:Mratio}	\\ \hline															
$	40	^2$&	42.0	&$	2.25	e{	-2	}$&	--	&$	2.24	e{	-2	}$&	--	&$	5.25	e{	-2	}$&	--	\\			
$	80	^2$&	338.1	&$	2.94	e{	-4	}$&	6.26	&$	2.77	e{	-4	}$&	6.34	&$	6.80	e{	-4	}$&	6.27	\\	
$	160	^2$&	2704.0	&$	2.81	e{	-6	}$&	6.71	&$	2.75	e{	-6	}$&	6.66	&$	1.05	e{	-5	}$&	6.02	\\			
$	320	^2$&	21730.0	&$	3.53	e{	-8	}$&	6.31	&$	3.50	e{	-8	}$&	6.30	&$	1.67	e{	-7	}$&	5.98	\\			
\hline\hline																									
\multicolumn{8}{|c|}{{\bf 2D RIDG: ($\mdeg=5$, $\nu=0.75$)}}	\\ \hline
{\bf mesh}		&	{\bf$T_r $ (s)}	&	$L^1$ {\bf error}	& \cref{eqn:Mratio}	&	$L^2$ {\bf error}	& \cref{eqn:Mratio}	
	&	$L^\infty$ {\bf error}	& \cref{eqn:Mratio}	\\ \hline																					
$	40	^2$&	4.1	&$	5.76	e{	-3	}$&	--	&$	5.86	e{	-3	}$&	--	&$	3.30	e{	-2	}$&	--	\\			
$	80	^2$&	32.9	&$	1.62	e{	-4	}$&	5.15	&$	1.54	e{	-4	}$&	5.25	&$	6.30	e{	-4	}$&	5.71	\\			
$	160	^2$&	254.5	&$	2.25	e{	-6	}$&	6.18	&$	2.16	e{	-6	}$&	6.15	&$	9.18	e{	-6	}$&	6.10	\\			
$	320	^2$&	2035.1	&$	3.04	e{	-8	}$&	6.21	&$	3.00	e{	-8	}$&	6.17	&$	1.51	e{	-7	}$&	5.93	\\			
\hline																									
\end{tabular}
\caption{
Convergence and runtime study for the 2D LIDG and RIDG  methods with $\mdeg=3$ and $\mdeg=5$
on 2D advection equation \cref{eqn:adv2d} with initial condition \cref{cc2D}.  Shown for various mesh sizes are the runtimes 
($T_r$) measured in seconds, 
the relative errors in $L^1$, $L^2$, and $L^{\infty}$, as well as the estimated convergence rates according to formula \cref{eqn:Mratio}.  We see that for any fixed number of elements, the RIDG method has a shorter runtime (i.e., computational cost) and slightly smaller error.
\label{table:RIvLI2D}}
\end{table}

\subsection{2D convergence tests}
We consider the 2D advection equation \cref{eqn:adv2d} with $u_x=u_y=1$, $\Omega=[-1,1]^2$, double periodic BCs,
and initial condition: 
\begin{equation}
\label{cc2D}
q(t=0,x,y) = \sin(16 \pi x)\sin(16 \pi y).
\end{equation}
We run the code \cite{code:ridg-code} to $t=2$ with $\mdeg=3$ (LIDG: $\nu=0.05$, RIDG: $\nu=0.75$) and $\mdeg=5$ (LIDG: $\nu=0.03$, RIDG: $\nu=0.75$) and compare runtimes and errors;  the results are shown in \cref{table:RIvLI2D}. 
We again see that both methods exhibit the expected convergence rates in 
$L^1$, $L^2$, and $L^\infty$.
For a fixed number of elements, usage of the RIDG method leads to slightly smaller errors.   We also notice that for a fixed number of elements the experiment runtime for the RIDG method is shorter than that of the LIDG method -- this is due to the increase in the maximum linearly stable CFL number from LIDG to RIDG.


\begin{table}[!t]																									
\centering																									
\begin{tabular}{|r|r|c|c|c|c|c|c|}																									
\hline																									
\multicolumn{8}{|c|}{{\bf 3D LIDG: ($\mdeg=3$, $\nu=0.03$)}}	\\ \hline																					
{\bf mesh}		&	{\bf $T_r$ (s)}	&	$L^1$ {\bf error}	& \cref{eqn:Mratio}	&	$L^2$ {\bf error}	& \cref{eqn:Mratio}							
	&	$L^\infty$ {\bf error}	& \cref{eqn:Mratio}	\\ \hline																					
$	20	^2$&	123.9	&$	1.21	e{	-3	}$&	--	&$	1.20	e{	-3	}$&	--	&$	6.16	e{	-3	}$&	--	\\			
$	40	^3$&	2027.4	&$	6.82	e{	-5	}$&	4.15	&$	6.95	e{	-5	}$&	4.11	&$	3.92	e{	-4	}$&	3.97	\\			
$	80	^3$&	32632.4	&$	4.21	e{	-6	}$&	4.02	&$	4.31	e{	-6	}$&	4.01	&$	2.50	e{	-5	}$&	3.97	\\			
\hline\hline																									
\multicolumn{8}{|c|}{{\bf 3D RIDG: ($\mdeg=3$, $\nu=0.6$)}}	\\ \hline																								
{\bf mesh}		&	{\bf $T_r$ (s)}	&	$L^1$ {\bf error}	& \cref{eqn:Mratio}	&	$L^2$ {\bf error}	& \cref{eqn:Mratio}																
	&	$L^\infty$ {\bf error}	& \cref{eqn:Mratio}	\\ \hline																					
$	20	^3$&	28.5	&$	9.24	e{	-4	}$&	--	&$	9.86	e{	-4	}$&	--	&$	5.02	e{	-3	}$&	--	\\			
$	40	^3$&	457.8	&$	5.85	e{	-5	}$&	3.98	&$	6.21	e{	-5	}$&	3.99	&$	3.15	e{	-4	}$&	3.99	\\			
$	80	^3$&	7299.2	&$	3.68	e{	-6	}$&	3.99	&$	3.89	e{	-6	}$&  6.20	&$	1.96	e{	-5	}$&	4.01	\\			
\hline\hline																									
\multicolumn{8}{|c|}{{\bf 3D LIDG: ($\mdeg=5$, $\nu=0.025$)}}	\\ \hline																								
{\bf mesh}		&	{\bf $T_r$ (s)}	&	$L^1$ {\bf error}	& \cref{eqn:Mratio}	&	$L^2$ {\bf error}	& \cref{eqn:Mratio}																
	&	$L^\infty$ {\bf error}	& \cref{eqn:Mratio}	\\ \hline																					
$	20	^3$&	318.8	&$	1.41	e{	-5	}$&	--	&$	1.34	e{	-5	}$&	--	&$	9.19	e{	-5	}$&	--	\\			
$	40	^3$&	4960.1	&$	1.91	e{	-7	}$&	6.21	&$	1.79	e{	-7	}$&	6.23	&$	1.53	e{	-6	}$&	5.90	\\			
$	80	^3$&	77582.5	&$	2.65	e{	-9	}$&	6.17	&$	2.50	e{	-9	}$&	6.16	&$	2.48	e{	-8	}$&	5.95	\\			
\hline\hline																									
\multicolumn{8}{|c|}{{\bf 3D RIDG: ($\mdeg=5$, $\nu=0.6$)}}	\\ \hline																						
{\bf mesh}		&	{\bf $T_r$ (s)}	&	$L^1$ {\bf error}	& \cref{eqn:Mratio}	&	$L^2$ {\bf error}	& \cref{eqn:Mratio}							
	&	$L^\infty$ {\bf error}	& \cref{eqn:Mratio}	\\ \hline																					
$	20	^3$&	61.6	&$	1.01	e{	-5	}$&	--	&$	9.77	e{	-6	}$&	--	&$	6.91	e{	-5	}$&	--	\\			
$	40	^3$&	960.2	&$	1.37	e{	-7	}$&6.20	&$	1.36	e{	-7	}$&	6.16	&$	1.08	e{	-6	}$&	6.00	\\			
$	80	^3$&	15510.1	&$	1.81	e{	-9	}$&	6.24	&$	1.88	e{	-9	}$&	6.18	&$	1.72	e{	-8	}$&	5.98	\\			
\hline																									
\end{tabular}																									
\caption{Convergence and runtime study for the 3D LIDG and RIDG  methods with $\mdeg=3$ and $\mdeg=5$														
on 3D advection equation \cref{eqn:adv3d} with initial condition \cref{cc3D}.  Shown for various mesh sizes are the runtimes									($T_r$) measured in seconds,																the relative errors in $L^1$, $L^2$, and $L^{\infty}$, as well as the estimated convergence rates according to formula \cref{eqn:Mratio}.  We see that for any fixed number of elements, the RIDG method has a shorter runtime (i.e., computational cost) and slightly smaller error.
\label{table:RIvLI3D}																									
}																									
\end{table}%

\subsection{3D convergence tests}
We consider the 3D advection equation \cref{eqn:adv3d} with $u_x=u_y=u_z=1$, $\Omega=[-1,1]^3$, triple periodic BCs,
and initial condition: 
\begin{equation}
\label{cc3D}
q(t=0,x,y,z) = \sin(2 \pi x)\sin(2 \pi y)\sin(2 \pi z).
\end{equation}
We run the code \cite{code:ridg-code} to $t=2$ with $\mdeg=3$ (LIDG: $\nu=0.03$, RIDG: $\nu=0.6$) and $\mdeg=5$ (LIDG: $\nu=0.025$, RIDG: $\nu=0.6$) and compare the error properties of the solution produced by the LIDG and RIDG methods;  the results are shown in \cref{table:RIvLI3D}.
We again see that both methods exhibit the expected convergence rates in 
$L^1$, $L^2$, and $L^\infty$.
As in the one and two-dimensional settings, the RIDG method exhibits better error and runtime properties than the LIDG method.


\section{Nonlinear RIDG}
\label{sec:burgers}

We show in this section how to extend the regionally-implicit discontinous Galerkin scheme (RIDG) to nonlinear problems. We show computational comparisons of the proposed RIDG scheme to a standard Runge-Kutta discontinuous Galerkin (RKDG) scheme on the 1D and 2D Burgers equation.

\subsection{Burgers equation in 1D}
We  consider the nonlinear inviscid Burgers equation in 1D:
\begin{equation}
\label{eqn:burgers1D}
q_{,t} + \frac{1}{2}\paren{q^2}_{,x}  = 0,
\end{equation}
where $(t,x) \in \left[ 0, T \right] \times \left [0,2\pi \right ]$ and
periodic boundary conditions are assumed. $T$ is chosen as
some time before shock formation occurs in the exact solution. 
 The initial conditions are taken to be
\begin{equation}
\label{eqn:burgers1D_IC}
q(t=0,x) = 1 - \cos(x).
\end{equation}

Unlike in the linear advection case, nonlinear conservation laws will require us to solve nonlinear algebraic equations in each of the regions depicted in \cref{fig:RIDG_1D}.  These nonlinear algebraic equations can be written in terms
of a nonlinear residual defined on each region:
\begin{equation}
\label{eqn:ridg1d_residual_total}
\vec {\mathcal R} = \vectthree{\vec {\mathcal R_1}}{\vec {\mathcal R_2}}{\vec {\mathcal R_3}},
\end{equation}
where
\begin{align}
\label{eqn:ridg1d_residual_1}
\begin{split}
\vec {\mathcal R_1} = & \int  \psifutr\psifutr^T \Wm  d\vec{\mathcal S}_{\tau}
- \int  \psipast\psifutr^T \Wmpast  d\vec{\mathcal S}_{\tau} \\
& + \nu_x\int \psieast \tilde f\paren{  \psieast^T \Wm ,  \psiwest^T \Wo } d{\vec{\mathcal S}}_{\xi} \\
& - \nu_x\int \psiwest f\paren{ \psiwest^T \Wm } d{\vec{\mathcal S}}_{\xi} \\
& - \iint  \psitau  \vec{\Psi}^T \Wo  d\vec{\mathcal S}  -\nu_x\iint  \psixii  f\paren{\vec{\Psi}^T \Wo}  d\vec{\mathcal S}, 
\end{split} \\
\label{eqn:ridg1d_residual_2}
\begin{split}
\vec {\mathcal R_2} = &\int  \psifutr\psifutr^T \Wo  d\vec{\mathcal S}_{\tau}
- \int  \psipast\psifutr^T \Wopast  d\vec{\mathcal S}_{\tau} \\
& + \nu_x\int \psieast \tilde f\paren{  \psieast^T \Wo ,  \psiwest^T \Wp } d{\vec{\mathcal S}}_{\xi} \\
& - \nu_x\int \psiwest \tilde f\paren{  \psieast^T \Wm ,  \psiwest^T \Wo } d{\vec{\mathcal S}}_{\xi} \\
& - \iint  \psitau \vec{\Psi}^T \Wo  d\vec{\mathcal S}  -\nu_x \iint  \psixii  f\paren{\vec{\Psi}^T \Wo}  d\vec{\mathcal S},
\end{split} \\
\label{eqn:ridg1d_residual_3}
\begin{split}
\vec {\mathcal R_3}= & \int  \psifutr\psifutr^T \Wp  d\vec{\mathcal S}_{\tau}
- \int  \psipast\psifutr^T \Wppast  d\vec{\mathcal S}_{\tau} \\
&+ \nu_x \int \psieast f\paren{ \psieast^T \Wp } d{\vec{\mathcal S}}_{\xi} \\
& - \nu_x\int \psiwest \tilde f\paren{  \psieast^T \Wo ,  \psiwest^T \Wp } d{\vec{\mathcal S}}_{\xi} \\
 &- \iint  \psitau  \vec{\Psi}^T \Wp  d\vec{\mathcal S}  -\nu_x \iint  \psixii  f\paren{\vec{\Psi}^T \Wp}  d\vec{\mathcal S}.
 \end{split}
\end{align}

For such general nonlinear problems, we use the Rusanov \cite{article:Ru61} numerical flux in lieu of the upwinded fluxes  in the prediction step (seen in \cref{fig:RIDG_1D}) and again for the time-averaged fluxes in the correction step.  For a scalar conservation law with flux function $f(q)$ and flux Jacobian $f'(q)$, the Rusanov flux is 
\begin{equation}
  \label{eqn:rusanov} 
 \tilde f(q_{\ell}, q_r) = \frac{1}{2} \left( f(q_{\ell}) + f(q_r) \right)- \frac{\lambda
 \left( q_{\ell}, q_r \right)}{2} \paren{ q_r - q_{\ell} },
 \end{equation}
where for scalar conservation laws:
 \begin{equation}
 \lambda\left( q_{\ell}, q_r \right) = \max \left\{\abs{f'(q_{\ell})},\abs{f'((q_{\ell}+q_r)/2)},\abs{f'(q_r)}\right\}.
 \end{equation}
  For Burgers equation \cref{eqn:burgers1D}, the numerical flux \cref{eqn:rusanov} becomes 
\begin{equation}
 \tilde f(q_{\ell}, q_r) =  \frac{1}{4}q_{\ell}^2 + \frac{1}{4}q_r^2 - \frac{\max \left\{ \abs{q_{\ell}},\abs{q_r} \right\}}{2} \left(q_r - q_{\ell}\right).
 \end{equation}
 
The goal in each region in each time-step is to minimize residual \cref{eqn:ridg1d_residual_total} with 
respect to the unknown space-time Legendre coefficients of the approximate
solution. We accomplish this by utilizing a Newton iteration. When forming the Newton iteration Jacobian (not to be confused with the flux Jacobian of the hyperbolic conservation law), one must compute the Jacobian of \cref{eqn:ridg1d_residual_total} by differentiating with respect to each coefficient. That is, we must compute 
\begin{equation}
\label{eqn:ridg1d_jacobian_total}
\mat{\mathcal J} = %
\begin{pmatrix}[2]
 \frac{\partial \vec {\mathcal R_1}}{\partial \Wm} &
 \frac{\partial \vec {\mathcal R_1}}{\partial \Wo} & 
 0 \\
 \frac{\partial \vec {\mathcal R_2}}{\partial \Wm} &
 \frac{\partial \vec {\mathcal R_2}}{\partial \Wo} &
 \frac{\partial \vec {\mathcal R_2}}{\partial \Wp} \\
 0 &
 \frac{\partial \vec {\mathcal R_3}}{\partial \Wo} &
 \frac{\partial \vec {\mathcal R_3}}{\partial \Wp}
 \end{pmatrix},
\end{equation}
which is analogous to the coefficient matrix in the linear advection case
(e.g., see \cref{eqn:ridg1d_system}). 
When computing the entries in \cref{eqn:ridg1d_jacobian_total}, one must deal with
the fact that the wave speed that appears in the Rusanov flux \cref{eqn:rusanov} is not
a smooth function of the coefficients; in order to handle this
issue we impose in
the computation of Jacobian \cref{eqn:ridg1d_jacobian_total} the
following condition:  $\frac{\partial}{\partial \vect W} \lambda = 0$.  This assumption seems to work well in practice, as evidenced in the results below, though other assumptions may be considered in the future. The stopping criterion for the Newton iteration that seems to be most effective at producing efficient solutions is the following
\begin{itemize}
    \item Stop if the residual of the region's main cell (the cell for whom we are forming a prediction) is below a certain tolerance ($\text{TOL} = 10^{-4}$);
    \item Stop if a maximum number of iterations is reached ($N_{\text{iters}} = 3$).
\end{itemize}

Unlike in the linear case, for  nonlinear conservation laws we cannot completely precompute the prediction update. However, we are able to leverage a so-called quadrature-free implementation \cite{Atkins1998} to increase the efficiency of the quadrature part of the prediction step. We review this methodology here: when computing the space-time quadrature, the most computationally expensive pieces are associated with the term in the residual and Jacobian where we integrate in space-time over the cell. 

For example, one such term for Burgers equation is:
\begin{equation} 
\begin{split}
\vec {\mathcal R_1}_{|\text{volume}} &= -\nu_x\iint  \vec{\Psi}_{,\xi}  \, f\paren{\vec{\Psi}^T \Wm} \, d\vec{\mathcal S} \\
& = - \nu_x\iint  \vec{\Psi}_{,\xi} \, \frac{1}{2}\paren{\vec{\Psi}^T \Wm}^2 \, d\vec{\mathcal S}.
\end{split}
\end{equation}
The contribution to the Jacobian matrix from this term is 
\begin{equation}
 \label{expensive}
\frac{ \partial {\mathcal R_1}_{|\text{volume}}}{ \partial \Wm }  = -\nu_x \iint \vec\Psi_{,\xi} \, \vec{\Psi}^T \Wm \, \vec\Psi^T  \, d\vec S.
\end{equation}
Each entry of this matrix has the form:
\begin{equation}
\label{eq:quadfree_expressions}
\left [ \frac{ \partial {\mathcal R_1}_{|\text{volume}}}{ \partial \Wm } \right ]_{ab}  
= \iint \psi^{(a)}_{,\xi} \, \sum\limits_{\ell = 1}^{\theta_T} \psi^{(\ell)} \, Q_{\ell} \psi^{(b)} \,  d\vec S     
= \paren{ \sum\limits_{\ell = 1}^{\theta_T}  \iint \psi^{(a)}_{,\xi} \psi^{(\ell)}  \psi^{(b)}   d\vec S }    Q_{\ell}.
\end{equation}
Notice that the expression within the parentheses found in \cref{eq:quadfree_expressions} can be precomputed using exact expressions (e.g., with a symbolic toolbox if the exact expressions are too laborious to derive by hand).  This allows us to forego an expensive quadrature routine in favor of an exact expansion of the coefficients $Q_\ell$ for forming the Newton iteration Jacobian.  In \cite{Shu1988} this idea was effectively expanded to certain types of non-polynomial flux functions, indicating that this idea can be generalized.  Thus we can avoid space-time quadrature of the volume integrals by integrating the expressions such as in \cref{eq:quadfree_expressions} analytically.


\begin{table}[!t]
\centering
\begin{tabular}{|r|r|r|c|c|c|c|c|c|}
\hline	
\multicolumn{9}{|c|}{{\bf 1D RKDG: ($\mdeg=3$, $\nu=0.1$)}}	\\ \hline
{\bf mesh}		& $N_T$ &	{\bf $T_r$ (s)}	&	$L^1$ {\bf error}	& \cref{eqn:Mratio}	&	$L^2$ {\bf error}	& \cref{eqn:Mratio}
	&	$L^\infty$ {\bf error}	& \cref{eqn:Mratio} 	\\ \hline
39	&	30	&	0.404	&	1.45E-07	&	-	&	2.39E-07	&	-	&	1.30E-06	&	-	\\
52	&	39	&	0.657	&	4.69E-08	&	3.91	&	7.69E-08	&	3.94	&	4.26E-07	&	3.89	\\
65	&	48	&	0.879	&	1.95E-08	&	3.94	&	3.18E-08	&	3.95	&	1.77E-07	&	3.93	\\
77	&	57	&	1.243	&	9.93E-09	&	3.98	&	1.63E-08	&	3.96	&	9.08E-08	&	3.95	\\
91	&	66	&	1.732	&	5.11E-09	&	3.98	&	8.39E-09	&	3.96	&	4.68E-08	&	3.97	\\
105	&	76	&	2.283	&	2.94E-09	&	3.86	&	4.76E-09	&	3.97	&	2.66E-08	&	3.96	\\
158	&	114	&	5.138	&	5.72E-10	&	4.01	&	9.36E-10	&	3.98	&	5.23E-09	&	3.98	\\
\hline\hline
\multicolumn{9}{|c|}{{\bf 1D RIDG: ($\mdeg=3$, $\nu=0.9$)}}	\\ \hline
{\bf mesh}	& $N_T$	&	{\bf $T_r$ (s)}	&	$L^1$ {\bf error}	& \cref{eqn:Mratio}	&	$L^2$ {\bf error}	& \cref{eqn:Mratio} &	$L^\infty$ {\bf error}	& \cref{eqn:Mratio} \\ \hline
39	&	3	&	0.125	&	1.47E-07	&	-	&	2.35E-07	&	-	&	1.48E-06	&	-	\\
52	&	4	&	0.217	&	4.70E-08	&	3.97	&	7.55E-08	&	3.94	&	4.85E-07	&	3.88	\\
65	&	5	&	0.330	&	1.93E-08	&	4.00	&	3.12E-08	&	3.96	&	2.01E-07	&	3.95	\\
77	&	6	&	0.501	&	9.69E-09	&	4.05	&	1.61E-08	&	3.93	&	1.06E-07	&	3.79	\\
91	&	7	&	0.578	&	4.95E-09	&	4.03	&	8.24E-09	&	4.00	&	5.65E-08	&	3.75	\\
105	&	8	&	0.810	&	2.82E-09	&	3.93	&	4.69E-09	&	3.94	&	3.24E-08	&	3.89	\\
158	&	12	&	2.136	&	5.54E-10	&	3.98	&	9.26E-10	&	3.97	&	6.50E-09	&	3.93	\\
\hline\hline
\multicolumn{9}{|c|}{{\bf 1D RIDG: ($\mdeg=5$, $\nu=0.9$)}}	\\ \hline
{\bf mesh}	& $N_T$	&	{\bf $T_r$ (s)}	&	$L^1$ {\bf error}	& \cref{eqn:Mratio}	&	$L^2$ {\bf error}	& \cref{eqn:Mratio} &	$L^\infty$ {\bf error}	& \cref{eqn:Mratio} \\ \hline
13	&	1	&	0.085	&	4.03E-08	&	-	&	6.79E-08	&	-	&	4.98E-07	&	-	\\
26	&	2	&	0.303	&	6.90E-10	&	5.87	&	1.20E-09	&	5.82	&	9.38E-09	&	5.73	\\
39	&	3	&	0.653	&	6.73E-11	&	5.74	&	1.22E-10	&	5.64	&	1.35E-09	&	4.79	\\
53	&	4	&	0.958	&	1.03E-11	&	6.11	&	1.75E-11	&	6.31	&	1.81E-10	&	6.54	\\
66	&	5	&	1.437	&	2.68E-12	&	6.15	&	4.77E-12	&	5.93	&	5.33E-11	&	5.57	\\\hline\hline
\multicolumn{9}{|c|}{{\bf 1D RIDG: ($\mdeg=7$, $\nu=0.9$)}}	\\ \hline
{\bf mesh}	& $N_T$	&	{\bf $T_r$ (s)}	&	$L^1$ {\bf error}	& \cref{eqn:Mratio}	&	$L^2$ {\bf error}	& \cref{eqn:Mratio} &	$L^\infty$ {\bf error}	& \cref{eqn:Mratio} \\ \hline
3	&	1	&	0.046	&	1.31E-05	&	-	&	2.43E-05	&	-	&	1.06E-04	&	-	\\
8	&	1	&	0.192	&	9.53E-09	&	7.37	&	1.48E-08	&	7.55	&	1.22E-07	&	6.90	\\ \hline
\end{tabular}  
    \caption{Convergence and runtime study for the 1D RKDG $\mdeg = 3$  and the RIDG methods for $\mdeg = 3, 5, 7$ on the inviscid Burgers equation \cref{eqn:burgers1D} with initial conditions \cref{eqn:burgers1D_IC}.  Shown for various mesh sizes are the number of time steps taken, $N_T$, the runtimes in seconds, $T_r$, and the relative errors as measured in the $L^1,L^2,$ and $L^\infty$ norms. }
    \label{fig:rkdg_compare}
\end{table}

Using this approach, we compared the nonlinear RIDG methods of various
orders to the $4^{\text{th}}$ order RKDG method as discussed in \cite{Shu1988}.
Shown in  \cref{fig:rkdg_compare} are the computed errors
and runtimes for the RKDG method  with $\mdeg = 3$ (i.e., the fourth-order RKDG
method) and the RIDG method $\mdeg = 3, 5, 7$; both methods were applied
to Burgers equation \cref{eqn:burgers1D} with initial conditions \cref{eqn:burgers1D_IC}, and run out to time $T=0.4$ (i.e, before shockwaves
form). We see that all methods exhibit the expected convergence rates in the $L^1$, $L^2$, and $L^\infty$ norms.  We note the following:
\begin{itemize}
    \item For any fixed error that we consider for the RKDG method, the RIDG method of $\mdeg=3$ can obtain a solution of similar accuracy about 2.5 to 3 times faster.  Furthermore, the RIDG method of $\mdeg=5$ can obtain a solution of similar accuracy 15 times faster.
   \item For any fixed error that we consider for the RKDG method, the RIDG method of $\mdeg=3$ can obtain a solution of similar accuracy while taking an order of magnitude fewer timesteps. 
   Furthermore, the RIDG method of $\mdeg=5$ can obtain a solution of similar accuracy while taking almost two orders of magnitude fewer timesteps. 
\end{itemize}

\subsection{Burgers equation in 2D}
Now we consider the nonlinear inviscid Burgers equation in 2D:
\begin{equation}
\label{eqn:burgers2D}
q_{,t} + \paren{\frac{1}{2}q^2}_{,x} + \paren{\frac{1}{2}q^2}_{,y}   =0,
\end{equation}
where $(t,x,y) \in \left [ 0, T \right] \times \left [0,2\pi \right]^2$. $T$ is some time before the shock forms in the solution.  We consider the initial conditions 
\begin{equation}
 \label{eqn:burgers2D_IC}
q(t=0,x) = \frac{1}{4}\bigl(1 - \cos(x)\bigr)\bigl(1-\cos(y)\bigr).
\end{equation}


\begin{table}[!t]
\centering
\begin{tabular}{|r|r|r|c|c|c|c|c|c|}
\hline	
\multicolumn{9}{|c|}{{\bf 2D RKDG: ($\mdeg=3$, $\nu=0.005$)}}	\\ \hline
{\bf mesh}		& $N_T$ &	{\bf $T_r$ (s)}	&	$L^1$ {\bf error}	& \cref{eqn:Mratio}	&	$L^2$ {\bf error}	& \cref{eqn:Mratio}
	&	$L^\infty$ {\bf error}	& \cref{eqn:Mratio}	\\ \hline
$	11	^2$	&	18	&	2.793	&	3.03e-05	&	--	&	4.44e-05	&	--	&	3.73e-04	&	--	\\
$	22	^2$	&	33	&	18.166	&	1.93e-06	&	3.97	&	2.97e-06	&	3.90	&	2.39e-05	&	3.96	\\
$	33	^2$	&	49	&	61.442	&	4.00e-07	&	3.88	&	6.11e-07	&	3.90	&	5.40e-06	&	3.67	\\
$	44	^2$	&	64	&	144.332	&	1.31e-07	&	3.88	&	1.98e-07	&	3.92	&	1.72e-06	&	3.98	\\
$	55	^2$	&	80	&	284.119	&	5.50e-08	&	3.89	&	8.21e-08	&	3.94	&	7.01e-07	&	4.01	\\
$	66	^2$	&	95	&	485.266	&	2.70e-08	&	3.91	&	3.99e-08	&	3.95	&	3.47e-07	&	3.85	\\
$	122	^2$	&	174	&	3075.613	&	2.37e-09	&	3.96	&	3.50e-09	&	3.96	&	2.99e-08	&	3.99	\\
\hline\hline
\multicolumn{9}{|c|}{{\bf 2D RIDG: ($\mdeg=3$, $\nu=0.75$)}}	\\ \hline
{\bf mesh}	& $N_T$	&	{\bf $T_r$ (s)}	&	$L^1$ {\bf error}	& \cref{eqn:Mratio}	&	$L^2$ {\bf error}	& \cref{eqn:Mratio}
	&	$L^\infty$ {\bf error}	& \cref{eqn:Mratio}	\\ \hline
$	11	^2$	&	1	&	6.021	&	2.89e-05	&	--	&	4.26e-05	&	--	&	2.76e-04	&	--	\\
$	22	^2$	&	2	&	45.500	&	1.85e-06	&	3.97	&	2.89e-06	&	3.88	&	1.87e-05	&	3.89	\\
$	33	^2$	&	3	&	142.633	&	3.93e-07	&	3.83	&	6.03e-07	&	3.86	&	4.38e-06	&	3.58	\\
$	44	^2$	&	4	&	315.483	&	1.29e-07	&	3.86	&	1.96e-07	&	3.91	&	1.33e-06	&	4.15	\\
$	55	^2$	&	5	&	585.140	&	5.43e-08	&	3.89	&	8.17e-08	&	3.92	&	5.84e-07	&	3.68	\\
$	66	^2$	&	6	&	984.351	&	2.66e-08	&	3.91	&	3.98e-08	&	3.94	&	2.75e-07	&	4.14	\\
$	122	^2$	&	12	&	6499.030	&	2.35e-09	&	3.95	&	3.53e-09	&	3.94	&	2.68e-08	&	3.78	\\
\hline\hline
\multicolumn{9}{|c|}{{\bf 2D RIDG: ($\mdeg=5$, $\nu=0.75$)}}	\\ \hline
{\bf mesh}	& $N_T$	&	{\bf $T_r$ (s)}	&	$L^1$ {\bf error}	& \cref{eqn:Mratio}	&	$L^2$ {\bf error}	& \cref{eqn:Mratio}
	&	$L^\infty$ {\bf error}	& \cref{eqn:Mratio}	\\ \hline
$	11	^2$	&	1	&	158.842	&	1.38e-07	&	--	&	2.42e-07	&	--	&	1.82e-06	&	--	\\
$	22	^2$	&	2	&	1190.925	&	2.35e-09	&	5.87	&	4.57e-09	&	5.73	&	4.53e-08	&	5.33	\\
$	33	^2$	&	3	&	3750.531	&	2.29e-10	&	5.74	&	4.41e-10	&	5.77	&	6.28e-09	&	4.87	\\
\hline
\end{tabular}
\caption{Convergence and runtime study for the 2D RKDG method, the 2D RIDG method for $\mdeg = 3$, and the RIDG method for $\mdeg = 5$ on the inviscid Burgers equation \cref{eqn:burgers2D} with initial conditions \cref{eqn:burgers2D_IC}.  Shown for various mesh sizes are the number of time steps taken, $N_T$, the runtimes in seconds, $T_r$, and the relative errors as measured in the $L^1,L^2,$ and $L^\infty$ norms. }
\label{fig:compare_2d}
\end{table}

We again use the Rusanov numerical flux for both the space-time surface integrals in the prediction step and the time-averaged fluxes in the correction step.  In \cref{fig:compare_2d} we compare the performance of the RKDG method to that of the RIDG method for $\mdeg=3, 5$ for $T = 0.4$. We observe the following
\begin{itemize}
    \item For all mesh sizes, the RKDG and RIDG methods with $\mdeg=3$ have similar error. The RIDG $\mdeg=3$ solutions take about 2 times longer to obtain, yet are obtained in 16 times fewer time steps. 
    \item For a fixed error of $\mathcal O (10^{-9})$, RKDG method takes 2.6 times longer to run than the the RIDG method $\mdeg=5$. Furthermore, the RIDG method $\mdeg=5$ obtains solutions with almost two orders of magnitude fewer time steps. 
\end{itemize}
We conclude that with respect to serial code execution in 2D, the RIDG method $\mdeg = 3$ is not as efficient as the RKDG method, but the RIDG method $\mdeg = 5$ is more efficient than the RKDG method.  This demonstrates the fact that the RIDG methods do not experience an analogous {\it Butcher barrier}, which causes efficiency deterioration for Runge-Kutta methods of orders higher than $4$.  Furthermore, with respect to minimizing the number of time steps (such as in the context of distributed memory programming), the RIDG method $\mdeg = 3$ require an order of magnitude fewer time steps for a fixed error than the RKDG method, while the RIDG method $\mdeg = 5$ requires almost two orders of magnitude fewer time steps for a fixed error. 
 
\section{Conclusions}
\label{sec:conclusions}
The purpose of this work was to develop a novel time-stepping method for high-order discontinuous Galerkin methods that has improved stability properties over traditional approaches (e.g., explicit SSP-RK and Lax-Wendroff).
The name we gave to this new approach is the {\it regionally-implicit} discontinuous Galerkin (RIDG) scheme, due to the fact that the prediction for a given cell is formed via an implicit method using information from small {\it regions} of cells around a cell, juxtaposed with the {\it local} predictor that forms a prediction for a given cell using only past information from that cell. More exactly, the RIDG method is comprised of a semi-localized version of a spacetime DG method, and a corrector step, which is an explicit method that uses the solution from the predictor step.  In this sense, the stencil of the RIDG schemes are slightly larger than similar explicit methods, and yet are able to take significantly larger time steps.

With this new scheme we achieved all of the following:
\begin{itemize}
\item Developed RIDG schemes for 1D, 2D, and 3D advection;
\item Demonstrated experimentally the correct convergence rates on 1D, 2D, and 3D advection examples;
\item Showed that the maximum linearly stable CFL number is bounded below by a constant that is independent of the polynomial order
(1D: 1.00, 2D: 0.75, 3D: 0.60);
\item Developed RIDG schemes for 1D and 2D nonlinear scalar equations;
\item Demonstrated experimentally the correct convergence rates on 1D and 2D nonlinear examples.
\item Showed that RIDG has larger maximum CFL numbers than explicit SSP-RKDG and Lax-Wendroff DG;
\item Demonstrated the efficiency of the RIDG schemes for nonlinear problems in 1D and 2D. Namely, we showed that the RIDG methods become more efficient as you increase the method order, as opposed to RKDG methods, whose efficiency deteriorates as you move beyond fourth-order accuracy due to the Butcher barrier.
\end{itemize}
All of the methods described in this work were written in a {\sc matlab} code that can be freely downloaded
\cite{code:ridg-code}.


There are many directions for future work for this class of methods, including
\begin{itemize}
    \item Exploring different methods for finding solutions to the nonlinear rootfinding problem that forms the predictions in each time step. This includes the possibility of constrained optimization so that the predictions fit some desired criterion such as maintaining positivity. 
    \item Extending the RIDG method to systems of equations while maintaining efficiency. A crucial development will be to extend existing limiter technology, including non-oscillatory limiters and positivity-preserving limiters, to the case of the RIDG scheme.  Since the RIDG method takes orders of magnitude fewer time steps when compared to SSP-RKDG and Lax-Wendroff DG, limiters that are used once or twice per time-step have a reduced effect on the overall runtime of the scheme.
    \item Implementing domain decomposition schemes and demonstrating efficient 
    many-core scaling for RIDG.  The RIDG stencil is small (almost like nearest neighbors), and would need to communicate only twice per time-step.  RKDG, Lax-Wendroff DG, and ADER-DG methods are known to be efficient on many-core systems \cite{Dumbser2018}, and so the RIDG method is similar enough to these methods that we expect similar results. However, since the RIDG method takes orders of magnitude fewer time steps when compared to these other methods, communication costs are minimized and so even greater efficiency might be achieved.  
\end{itemize}

\section*{Acknowledgments}
We would like to thank the anonymous referees for 
their thoughtful comments and suggestions that helped to improve this paper.  
This research was partially funded by NSF Grant DMS--1620128.

\bibliographystyle{siam}

\end{document}